%% file: main.tex
\documentclass[peerreview,letterpaper,12pt]{IEEEtran} 

\usepackage{setspace}
\ifCLASSINFOpdf
\else
\fi
\usepackage{float} 
\usepackage{graphicx}

\usepackage{subcaption}

\usepackage{datetime} 

\usepackage{multicol}
\usepackage{algorithm}
\usepackage{algpseudocode}
\usepackage{amsmath}
\usepackage{amssymb}
\usepackage{mathrsfs}
\usepackage{amsthm}
\usepackage{bm}
\usepackage{cite}
\usepackage{enumerate}
\usepackage{diagbox}
\usepackage{color}
\usepackage{comment}
\usepackage{epstopdf}

\usepackage{empheq} 
\usepackage{enumitem} 

\usepackage{tikz} 
\usetikzlibrary{calc} 

\usepackage{hyperref} 

\newcommand\ie{i.e.,~}

\newcommand{\cP}{\mathcal{P}}

\newcommand{\cH}{\mathcal{H}}

\newcommand{\bR}{\boldsymbol R}

\def\bR{\mathbb{R}} 

\newtheorem{theorem}{Theorem}
\newtheorem{definition}{Definition} 

\newtheorem{corollary}{Corollary}
\newtheorem{lemma}{Lemma}

\input{marc-macros.tex} 

\begin{document}
\title{Optimal Trajectories of a UAV Base Station Using Hamilton-Jacobi Equations}

\author{\IEEEauthorblockN{Marceau Coupechoux\IEEEauthorrefmark{1}, J\'er\^ome Darbon \IEEEauthorrefmark{2}, Jean-Marc K\'elif\IEEEauthorrefmark{3}, and Marc Sigelle\IEEEauthorrefmark{4}}
\IEEEauthorblockA{\IEEEauthorrefmark{1}LTCI, Telecom Paris, Institut Polytechnique de Paris, France, \IEEEauthorrefmark{2}Brown University, US, 
\IEEEauthorrefmark{3}Orange Labs, France, \IEEEauthorrefmark{4}On leave from Telecom Paris, France}
\IEEEauthorblockA{Email:  marceau.coupechoux@telecom-paristech.fr, jerome\_darbon@brown.edu, jeanmarc.kelif@orange.com,
marc.sigelle@gmail.com
\\
%
} 
}
\maketitle

\begin{abstract}
We consider the problem of optimizing the trajectory of an Unmanned Aerial Vehicle (UAV). Assuming a traffic intensity map of users to be served, the UAV must travel from a given initial location to a final position within a given duration and serves the traffic on its way. The problem consists in finding the optimal trajectory that minimizes a certain cost depending on the velocity and on the amount of served traffic. We formulate the problem using the framework of Lagrangian mechanics. We derive closed-form formulas for the optimal trajectory when the traffic intensity is quadratic (single-phase) using Hamilton-Jacobi equations. When the traffic intensity is bi-phase, i.e. made of two quadratics, we provide necessary conditions of optimality that allow us to propose a gradient-based algorithm and a new algorithm based on the linear control properties of the quadratic model. These two solutions are of very low complexity because they rely on fast convergence numerical schemes and closed form formulas. These two approaches return a trajectory satisfying the necessary conditions of optimality. At last, we propose a data processing procedure based on a modified K-means algorithm to derive a bi-phase model and an optimal trajectory simulation from real traffic data.     
\end{abstract}


\section{Introduction}
Unmanned Aerial Vehicles (UAV) are expected to play an increasing role in future wireless networks\footnote{J. Darbon is supported by NSF DMS-1820821. M. Coupechoux has performed his work at LINCS laboratory.}~\cite{zeng2019accessing}. UAVs may be deployed in an ad hoc manner when the traditional cellular infrastructure is missing. They can serve as relays to reach distant users outside the coverage of wireless networks. They also may be used to disseminate data to ground stations or collect information from sensors. 
In this paper, we address one of the envisioned use cases for UAV-aided wireless communications, which relates to cellular network offloading in highly crowded areas~\cite{zeng2019accessing}. More specifically, we focus on the path planning problem or trajectory optimization problem that consists in finding an optimal path for a UAV Base Station (BS) that minimizes a certain cost depending on the velocity and on the amount of served traffic. Our approach is based on the Lagrangian mechanics framework and the use of Hamilton-Jacobi (HJ) equations.

\subsection{Related Work}


UAV trajectory optimization for networks has been tackled maybe for the first time in~\cite{Pearre10}. The model consists in 
a UAV flying over a sensor network from which it has to collect some data. The UAV can learn from previous experience, which is not assumed in our study. 
The problem of optimally deploying UAV BSs to serve traffic demand has been addressed in the literature by considering static UAVs BSs or relays, see e.g.~\cite{sharma2016uav, yang2018energy}. The goal is to optimally position the UAV so as to maximize the data rate with ground stations or the number of served users. In these works, the notion of trajectory is either ignored or restricted to be circular or linear. 
In robotics and autonomous systems, trajectory optimization is known as {\it path planning}. In this field, there are classical methods like Cell Decomposition, Potential Field Method,  Probabilistic Road Map, or heuristic approaches, e.g. bio-inspired algorithms~\cite{Mac16}. Authors of \cite{chi2012civil} have capitalized on this literature and proposed a path planning algorithm for drone BSs based on A* algorithm. The main goal of these papers is to reach a destination while avoiding obstacles. In our work, we intend to minimize a certain cost function along the trajectory by controlling the velocity of the UAV.
This goal is studied in optimal control theory~\cite{Liberzon11} and is applied for example in aircraft trajectory planning~\cite{Delahaye14}. Most numerical methods in control theory can be classified in {\it direct} and {\it indirect} methods. In direct methods, the problem is transformed in a non linear programming problem using discretized time, locations and controls. 
Direct methods are heavily applied in a series of recent publications in the field of UAV-aided communications, see e.g.~\cite{Zeng16, Zeng17,wu2018joint,zeng2019energy}. Formulated problems are usually non-convex. The standard approach is hence to rely on Successive Convex Approximation (SCA), which iteratively minimizes a sequence of approximate convex functions. SCA is known to converge to a Karush-Kuhn-Tucker solution under mild conditions~\cite{marks1978general} but the quality of the solution may heavily depend on the initial guess. Here, simple heuristics or solutions to the Travelling Salesman Problem (TSP) or the Pickup-and-Deliver Problem (PDP) can be used for finding an initial trajectory~\cite{zhang2018uav}. With direct methods, because of the discretization, the differential equations and the constraints of the systems are satisfied only at discrete points. This can lead to less accurate solutions than indirect methods and the quality of the solution depends on the quantization step~\cite{vonStryk92}. Although every iteration of SCA has a polynomial time complexity, practical resolution time may dramatically increase with the quantization grid and the dimension of the problem. On the other hand, indirect approaches relies on considering the Hamilton-Jacobi Partial Differential Equation associated to the optimal control problem (see e.g., \cite{Barles16}, \cite{evans.10.book}[chp. 10]). Several recent methods have been proposed to solve HJ Partial Differential Equations (PDE) in high dimensions. These include  max-plus algebra methods \cite{akian2008max,mceneaney2006max}, dynamic programming and reinforcement learning \cite{alla2019efficient}, tensor decomposition techniques \cite{dolgov2019tensor}, sparse grids \cite{kang2017mitigating}, model order reduction \cite{alla2017error}, polynomial approximation \cite{kalise2019robust}, optimization methods \cite{chow.19.jcp,Darbon2016Algorithms,yegorov2017perspectives,chow.18.amsa} and neural networks \cite{kang.21.etal, Han2018Solving, darbonetal.20.rms, darbonMeng.20.jcp}. 
In this paper, we consider certain indirect methods that provide 
analytical solutions for certain classes of optimal control problem as we have shown in a preliminary study~\cite{coupechoux2019optimal}.

\subsection{Contributions}
Our contributions are the following:
\begin{itemize}
	\item {\it Problem Formulation:} To the best of our knowledge, this is the first time, after our preliminary study~\cite{coupechoux2019optimal}, that the UAV BS trajectory problem is formulated using the formalism of Lagrangian mechanics and solved using Hamilton-Jacobi equations. This approach provides closed-form equations when the potential is quadratic and thus very low complexity solutions compared to existing solutions in the literature.
    \item {\it Closed-form expression of the optimal trajectory with single phase traffic intensity:} When the traffic intensity map is made of a single hot spot or traffic hole, has a quadratic form ({\it single phase}), and is time-independent, closed form expressions for the optimal trajectory are derived. It follows a hyperbola for a hot spot and corresponds to a repulsor in mechanics. For a traffic hole, the trajectory is on an ellipse and corresponds to the case of an attractor in mechanics.  
    \item {\it Characterization of the optimal solution in multi-phase traffic intensity:} When the traffic map has several hot spots or traffic holes ({\it multi-phase}) whose regions are separated by interfaces and is time-independent, we derive necessary conditions to be fulfilled by the position and the instant at which the optimal trajectory crosses an interface (see Theorem~\ref{th:multiphase}). 
    \item {\it A gradient algorithm for bi-phase traffic:} An in-depth analysis of
    convexity vs. non-convexity
    issues allows us to derive a gradient algorithm to solve the bi-phase problem (Algorithm~\ref{alg:ms-grad-algo}). 
    This algorithm finds a stationary point for the cost function. This algorithm has a complexity $O(1)$ at every iteration, whereas iterations of the sequential convex optimization technique have polynomial time complexity.
    \item {\it A new algorithm for the bi-phase optimization problem:} A new algorithm, called the $B$-algorithm (Algorithm~\ref{alg:bcurve-dichotomy}), is proposed based on the linear control properties of the quadratic model. 
    This algorithm relies on a bisection scheme the complexity of which is proportional to the logarithm of the desired precision and closed form formulas.
    \item {\it A data processing procedure:} We propose a method to pre-process real measured traffic data in order to derive a bi-phase quadratic model. This procedure is based on smoothing steps followed by a modified K-means algorithm adapted to our quadratic model (Algorithm~\ref{alg:preproc}). In our numerical experiments, the optimal trajectory is computed in a region where real traffic data is available~\cite{Modeling15-Chen}.   
\end{itemize}
The paper is structured as follows. In Section
~\ref{sec:systemmodel}, we give the system model, its interpretation in terms of Lagrangian mechanics
formulate the problem and give preliminary results. Section~\ref{sec:quadcost} is devoted to the characterization of the optimal trajectories
for both the single- and bi-phase cases. Section~\ref{sec:algos} presents our algorithms, Section~\ref{sec:numericalexp} our data processing procedure and numerical experiments. Section~\ref{sec:conclusion} concludes the paper.  
\\
{\bf Notations:} 
The usual Euclidean scalar product between $x\in \bR^n$ and $y\in \bR$ is denoted by $x \cdot y$.
The Euclidean norm $\|x\|$ in $\mathbb{R}^n$ of $x\in\mathbb{R}^n$ is defined by $\|x\| := \sqrt{x \cdot x}$.
The set of matrices with $m$ rows, $n$ columns and real entries is denoted by $\mathcal{M}_{m,n}(\bR)$.
The transpose of the $A\in \mathcal{M}_{m,n}(\bR)$ is denoted by 
$\transp{A} \in \mathcal{M}_{n,m}(\bR)$.
We classically identify $\mathcal{M}_{m,1}(\bR)$ and $\mathcal{M}_{1,n}(\bR)$ as column vectors of $\bR^m$ and row vectors of $\bR^n$, respectively. 
Let $f:\mathbb{R}^n\times \mathbb{R}^m \to \mathbb{R}$ defined by $f(x,y)$ where $x = (x_1, \dots,x_n)\in \mathbb{R}^n$ and $y =(y_1,\dots,y_m)\in \mathbb{R}^m$. Let $a \in \mathbb{R}^n$ and $b \in \mathbb{R}^m$. We denote by $\frac{\partial f}{\partial x_i} (a,b)$ the partial derivative of $f$ with respect to the variable $x_i$ at $(a,b) \in \mathbb{R}^n \times \mathbb{R}^m$.
We also introduce the notations 
\\ 
$\nabla _x f(a,b) = ( \frac{\partial f}{\partial x_1}(a,b), \dots, \frac{\partial f}{\partial x_n}(a,b)) \in \mathbb{R}^n$ and $\nabla_y f(a,b) = ( \frac{\partial f}{\partial y_1}(a,b), \dots, \frac{\partial f}{\partial y_m}(a,b)) \in \mathbb{R}^m$.
\\
We also consider the following notation for partial Hessian matrices
\begin{equation}
\mathcal{M}_{m+n, m+n}(\bR) \ni \nabla^2 f(a,b)
=
\begin{pmatrix}
\nabla^2_{x,x} f(a,b) & \nabla^2_{x,y} f(a,b)\\
\nabla^2_{y,x} f(a,b) & \nabla^2_{y,y} f(a,b)
\end{pmatrix}
\label{eq:definition-hessian}
\end{equation}
where
$\nabla^2_{x,x} f(a,b) \in \mathcal{M}_{n,n}(\bR) = 
\begin{pmatrix}
\frac{\partial^2 f}{\partial x_1^2}(a,b)
&\dots &\frac{\partial^2 f}{\partial x_1 \partial x_n}(a,b)\\
\vdots & \dots & \vdots\\ 
\frac{\partial^2 f}{\partial x_n x_1}(a,b)
&\dots& \frac{\partial^2 f}{\partial x_n^2}(a,b) 
\end{pmatrix}$ and 
$\Id_n$ 
denotes the identity matrix of 
$\mathcal{M}_{n,n}(\bR)$.
We shall see that the value function $S:\bR \times \bR^2\times \bR\times \bR^2 \to \bR$ will play a fundamental role in this paper. We use the notations $(T_1, X_1,T_2, X_2)$ for $S$ and therefore the partial derivatives of $S$ at $(t_1, x_1, t_2, x_2) \in \bR\times \bR^2\times \bR\times \bR^2$ are denoted as follows: 
$
\dfrac{\partial S}{\partial T_1} (t_1, x_1, t_2, x_2)
$,
$
\nabla_{X_1} S (t_1, x_1, t_2, x_2)$,
$\dfrac{\partial S}{\partial T_2} (t_1, x_1, t_2, x_2)
$ 
and $
\nabla_{X_2} S (t_1, x_1, t_2, x_2)
$.

\section{System Model and
Lagrangian Mechanics 
Interpretation
\label{sec:systemmodel}}
\subsection{System Model}
\begin{figure}[t]
\begin{center}
\includegraphics[width=0.8\linewidth]{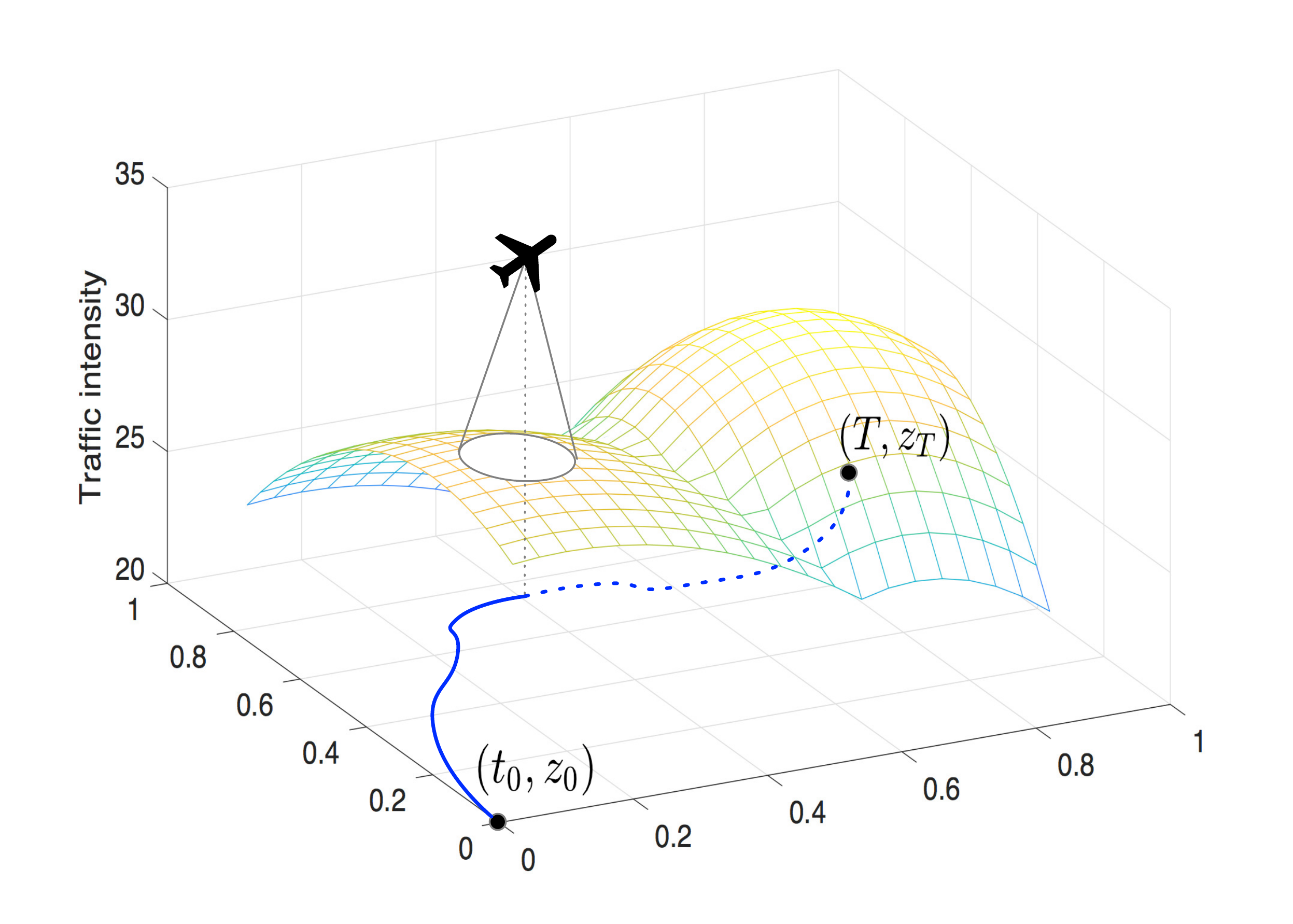}
\end{center}
\caption{\label{fig:systemmodel}A UAV Base Station travels from $z_0$ at $t_0$ to $z_T$ at $T$ and serves a user traffic characterized by its intensity.}
\end{figure}
We consider a network area characterized by a traffic density  at position $z$ and time $t$. We intend to control the trajectory and the velocity of a UAV base station, which is located in $z_0\triangleq z(t_0)$ at $t_0$ and shall reach a destination $z_T\triangleq z(T)$ at $T$ with the aim of minimizing a cost determined by the velocity and the traffic, defined hereafter by (\ref{eq:cost}). At $(t,z)$, we assume that the UAV BS is able to cover an area, from which it can serve users (see Figure~\ref{fig:systemmodel}). The velocity of the UAV BS induces an energy cost. In this model, we control the velocity vector $a$ of the UAV BS.
The general form of the cost function is as follows
\begin{equation} \label{eq:cost}
\mathcal{L}(t,z,a)=\frac{K}{2}||a||^2-u(t,z)
\end{equation}
where the first term is a cost related to the velocity of the vehicle ($K$ is a positive constant), and $\|\cdot\|$ denotes the usual Euclidean norm. The higher is the speed, the higher is the energy cost. The second term is a {\it user traffic intensity}, i.e., the amount of traffic served by the UAV BS at $(t,z)$. Note that a non-zero energy at null speed can be incorporated in the model by adding a constant. Without loss of generality, we assume that this constant is null.

Let $S(t_0, z_0,T,z_T)$ be the minimal total cost along any trajectory between $z_0$ at $t_0$ and $z_T$ at $T$ (also called {\it the action} in mechanics or {\it value function} in control theory).
Let us define $\Omega(t_0,T)$ as the space of absolutely continuous functions from $[t_0;T]$ to $\mathbb{R}^2$. Our problem can now be formulated as follows
\begin{eqnarray} \label{eq:problem}
S(t_0, z_0,T,z_T) \!\!\!&\!=\!\!\!\!&\!\!\!\!\! \min_{a\in \Omega(t_0,T) }\!\!\int_{t_0}^T\!\!\!\!
\mathcal{L}(s,z(s),a(s))ds \!\!+\!\! J(z(T))
\end{eqnarray}
where $\frac{d z}{d t}(t) = a(t)$, $z(t_0) = z_0$, and 
$J$ is the terminal cost defined by $J(z)=0$ if $z=z_T$ and $J(z)=+\infty$ otherwise.
For simplicity reasons, we assume the existence and uniqueness of the optimal control $a^*(t)$ in \eqref{eq:problem} and denote the associated optimal trajectory $z^*(t)$. In a traffic map symmetric with respect to $z_0$ and $z_T$, the reader can convince himself that the uniqueness is not guaranteed. 

\subsection{Preliminary Results From Lagrangian Mechanics}
We provide in this section important results from the Lagrangian mechanics for the convenience of the reader. 

\begin{definition}[Impulsion]
The \emph{impulsion} function is defined as 
\begin{equation}
p(t,z,a):=\nabla_a \mathcal{L}(t,z,a).
\end{equation}
\end{definition}
In the Newtonian classical framework that is used here (see \eqref{eq:cost}), the impulsion is the product of the particle mass by its velocity (hence the standard term ``impulsion").
\begin{definition}
The \emph{Hamiltonian} function is defined as
\begin{equation}
H(t,z,p):=\max_{a\in\mathbb{R}^2} p\cdot a-\mathcal{L}(t,z,a).
\end{equation}
\end{definition}

\begin{lemma}[Euler-Lagrange Equations] \label{lemma:EulerLagrange}
Along the optimal trajectory $z^*(t)$ that starts from $z_0$ at $t_0$ and ends at $z_T$ at $T$, we have
\begin{equation}\label{eq:eulerlagrange2}
\frac{d}{dt}\nabla_a \mathcal{L}(t,z^*(t),a^*(t))=\nabla_z \mathcal{L}(t,z^*(t),a^*(t))
\end{equation}
or equivalently
\begin{equation} \label{eq:eulerimpulsion}
\frac{dp}{dt}(t,z^*(t),a^*(t))=\nabla_z\mathcal{L}(t,z^*(t),a^*(t)).
\end{equation}
\end{lemma}

\begin{proof}
See Appendix~\ref{app:EulerLagrange}.
\end{proof}

The Euler-Lagrange equation is the first-order necessary condition for optimality and holds for every point on the optimal trajectory.  

\begin{lemma} \label{lemma:lagrangienhomogene}
If the Lagrangian $\mathcal{L}(t,z,a)$ is time-independent and $\alpha$-homogeneous in $z$ and $a$ for $\alpha>0$, i.e., $\mathcal{L}(\lambda z,\lambda a)=|\lambda|^{\alpha}\mathcal{L}(z,a)$ for all $\lambda\in \mathbb{R}$, $S$ given by \eqref{eq:problem} reads
\begin{equation} \label{eq:valuefunction}
S(t_0, z_0,T,z_T) = \frac{1}{\alpha}[z\cdot p]_{t_0}^T+J(z_T).
\end{equation}
\end{lemma}

\begin{proof}
See Appendix~\ref{app:lagrangienhomogene}.
\end{proof}

\begin{lemma}[Hamilton-Jacobi] \label{lemma:hamiltonjacobi} Along the optimal trajectory, we have for $t\in(t_0;T)$
using  previous notations 
\begin{equation}
\Hone{S}(t,z^*(t),T,z_T)
=H(t,z^*(t),-p^*(t)),
\label{Hamilton-Jacobi-backward-ms:eq}
\end{equation}
\begin{equation}
\Htwo{S} (t_0,z_0,t,z^*(t)) 
=-H(t,z^*(t),p^*(t)),
\label{Hamilton-Jacobi-forward-ms:eq}
\end{equation}
where 
\begin{equation} \label{eq:optimpulsion}
p^*(t)=
\nabla_a\mathcal{L}(t,z^*(t),a^*(t))
=
\Pone S(t_,z^*(t),T,z_T).
\end{equation}
\end{lemma}

\begin{proof}
See Appendix~\ref{app:hamiltonjacobi}.
\end{proof}

From now, we assume that the Lagrangian is time-independent, i.e., $\mathcal{L}(t,z,a)=\mathcal{L}(z,a)$, and is an even function in $a$, i.e., $\mathcal{L}(z,-a)=\mathcal{L}(z,a)$. A direct consequence is that $H$ is time-independent and is an even function in $p$, i.e., we write $H(t,z,p)=H(z,p)$ and $H(z,-p)=H(z,p)$. 

\section{Optimal Trajectory} \label{sec:quadcost}
In this section, we characterize the optimal trajectory when the traffic intensity is a quadratic form and also when it is made of two regions of quadratic form separated by an interface.
We call these two cases {\it single-phase} and {\it bi-phase} intensities respectively. Both cases satisfy our assumptions on the Lagrangian with $\alpha=2$.  

\subsection{Single-Phase Optimal Trajectory}

Assume that the traffic intensity is of the form $u(z)=\frac{1}{2}u_0||z||^2$. When $u_0>0$, this function models a traffic hole in $z=0$. When $u_0<0$, it models a traffic hot spot at $z=0$. We disregard the case $u_0=0$ because it corresponds to a constant traffic intensity that is not of interest in this paper. Thus the cost function has  the following form
\begin{equation} \label{eq:Lsinglephase}
\mathcal{L}(z,a)=\frac{1}{2}K||a||^2-\frac{1}{2}u_0||z||^2.
\end{equation}
Note that 
\begin{equation} \label{eq:impulsionquad}
p(z,a) = \nabla_a\mathcal{L}(z,a)=Ka.
\end{equation}

\subsubsection{Trajectory Equation} In the single phase case, we have a closed form expression of the trajectory. 

\begin{theorem} \label{th:quadraticfunction}
If $u_0<0$, the cost function is given by (\ref{eq:vu0neg}),
\begin{figure*}
   \begin{align} \label{eq:vu0neg}
S(t_0,z_0,T,z_T)=\frac{K\omega}{2\sinh \omega (T-t_0)}\left( (|z_0|^2+|z_T|^2)\cosh \omega (T-t_0)-2z_0\cdot z_T\right)+J(z_T) 
\end{align}
    \hrulefill
\end{figure*}   
the optimal trajectory is
\begin{equation} \label{eq:trajhotspot}
z^*(t)=\frac{z_T\sinh (\omega(t-t_0))+z_0\sinh (\omega(T-t))}{\sinh (\omega(T-t_0))}
\end{equation}
and the control is given by
\begin{equation} \label{eq:speedhotspot}
a^*(t)=\omega \frac{z_T\cosh(\omega(t-T))-z_0\cosh(\omega(T-t))}{\sinh(\omega(T-t_0))}
\end{equation}
where $\omega^2=-\frac{u_0}{K}$. \\
If $u_0>0$, the cost function is given by (\ref{eq:vu0pos}),
\begin{figure*}
   \begin{align} \label{eq:vu0pos}
S(t_0,z_0,T,z_T)=\frac{K\omega}{2\sin \omega (T-t_0)}\left( (|z_0|^2+|z_T|^2)\cos \omega (T-t_0)-2z_0\cdot z_T\right)+J(z_T)
\end{align}
    \hrulefill
\end{figure*}   
the optimal trajectory is
\begin{equation} \label{eq:ellipse}
z^*(t)=\frac{z_T\sin (\omega(t-t_0))+z_0\sin (\omega(T-t))}{\sin (\omega(T-t_0))}
\end{equation}
and the control is given by
\begin{equation}
a^*(t)=\omega \frac{z_T\cos (\omega(t-t_0))-z_0\cos (\omega(T-t))}{\sin (\omega(T-t_0))}
\end{equation}
where $\omega^2=\frac{u_0}{K}$. 
\end{theorem}

\begin{proof}
See Appendix~\ref{app:quadraticfunction}.
\end{proof}


\begin{corollary} \label{cor:quad}
If the user traffic intensity is of the form $u(t,z)=\frac{1}{2}u_0||z||^2+u_0z\cdot e +u_1$ with $u_1\in \mathbb{R}$ and $e\in \mathbb{R}^2$, then define $\tilde{z}=z+e$, $\tilde{z}_0=z_0+e$, $\tilde{z}_T=z_T+e$ and trajectories given in Theorem~\ref{th:quadraticfunction} are valid by replacing $z$, $z_0$, $z_T$ by $\tilde{z}$, $\tilde{z}_0$, $\tilde{z}_T$, respectively. The cost function becomes: 
$S(t_0,z_0,T,z_T) = \frac{1}{\alpha}[z\cdot p]_{t_0}^T+J(z_T)-u_1(T-t_0)$. 
\end{corollary}

\begin{corollary} \label{cor:barycentre}
If the user traffic intensity is of the form $u(t,z)=\sum_i u_i||z-z_i||^2$ with $\sum_i u_i\neq 0$, then $u(t,z)=(\sum_i u_i)||z-z_b||^2+\sum_i u_i||z_i-z_b||^2$ with $z_b=\frac{\sum_i u_iz_i}{\sum_i u_i}$. Define $\tilde{z}=z+z_b$, $\tilde{z}_0=z_0+z_b$, $\tilde{z}_T=z_T+z_b$, $\tilde{u}_0=\sum_i u_i$ and trajectories given in Theorem~\ref{th:quadraticfunction} are valid by replacing $z$, $z_0$, $z_T$, $u_0$ by $\tilde{z}$, $\tilde{z}_0$, $\tilde{z}_T$, $\tilde{u}_0$ respectively. 
\end{corollary}
The system is thus equivalent to the one assumed in Theorem~\ref{th:quadraticfunction} by changing the origin of the locations to the barycentre $z_b$ of the $z_i$.

\subsubsection{Traffic Hot Spot, Traffic Hole}
\label{subsub:hotspot-hole} 
\begin{figure}[t]
\begin{center}
\includegraphics[width=0.6\linewidth]{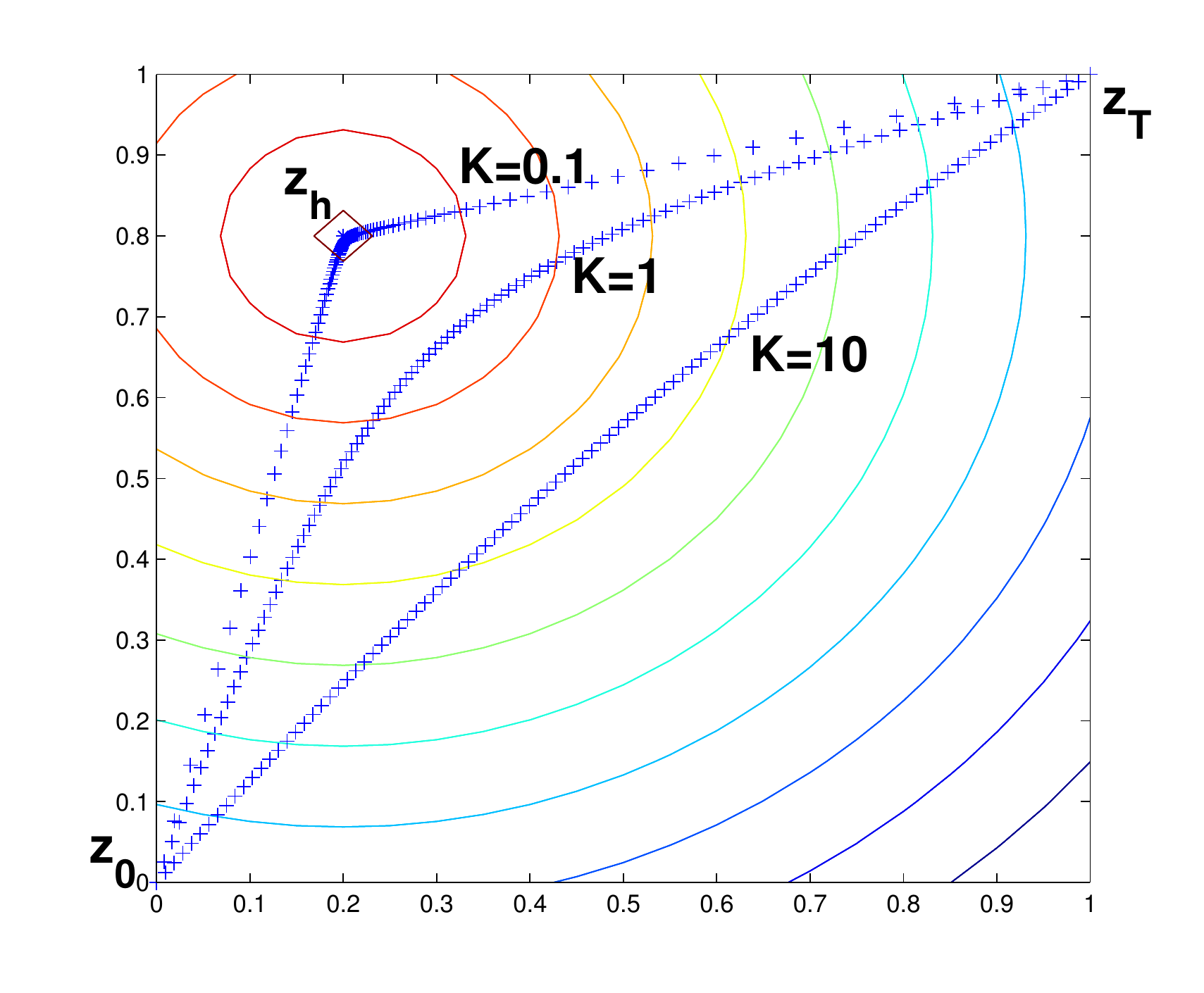}
\end{center}
\caption{\label{fig:hotspotu0negative}Traffic hot spot ($u_0<0$). Circles are iso-traffic levels.}
\end{figure}

\begin{figure*}[t!]
    \centering
    \begin{subfigure}[t]{0.5\textwidth}
        \centering
        \includegraphics[width=\linewidth]{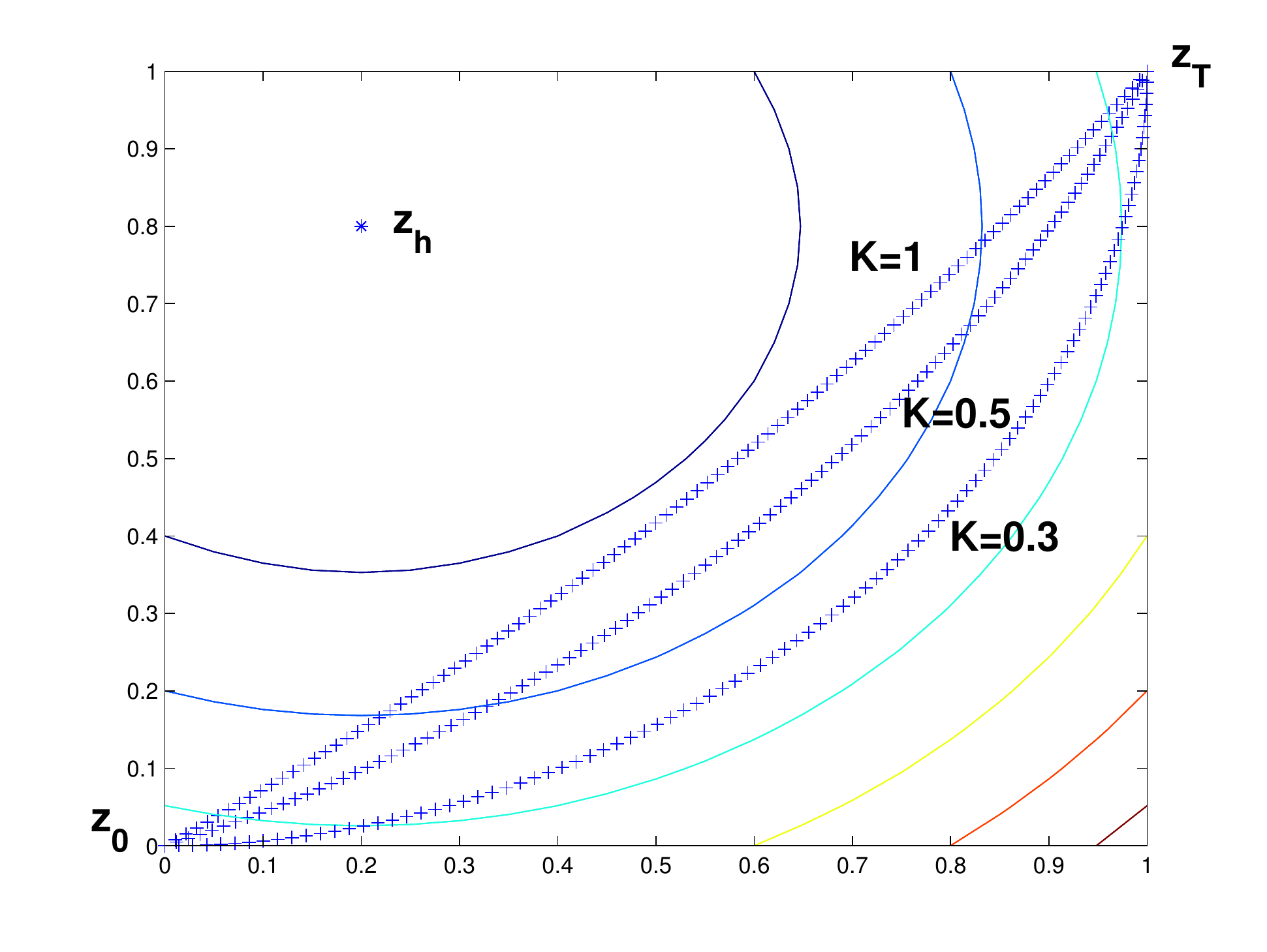}
        \caption{$T$ is smaller than the ellipse period.}
        \label{fig:u0positive2}
    \end{subfigure}%
    \begin{subfigure}[t]{0.5\textwidth}
        \centering
        \includegraphics[width=\linewidth]{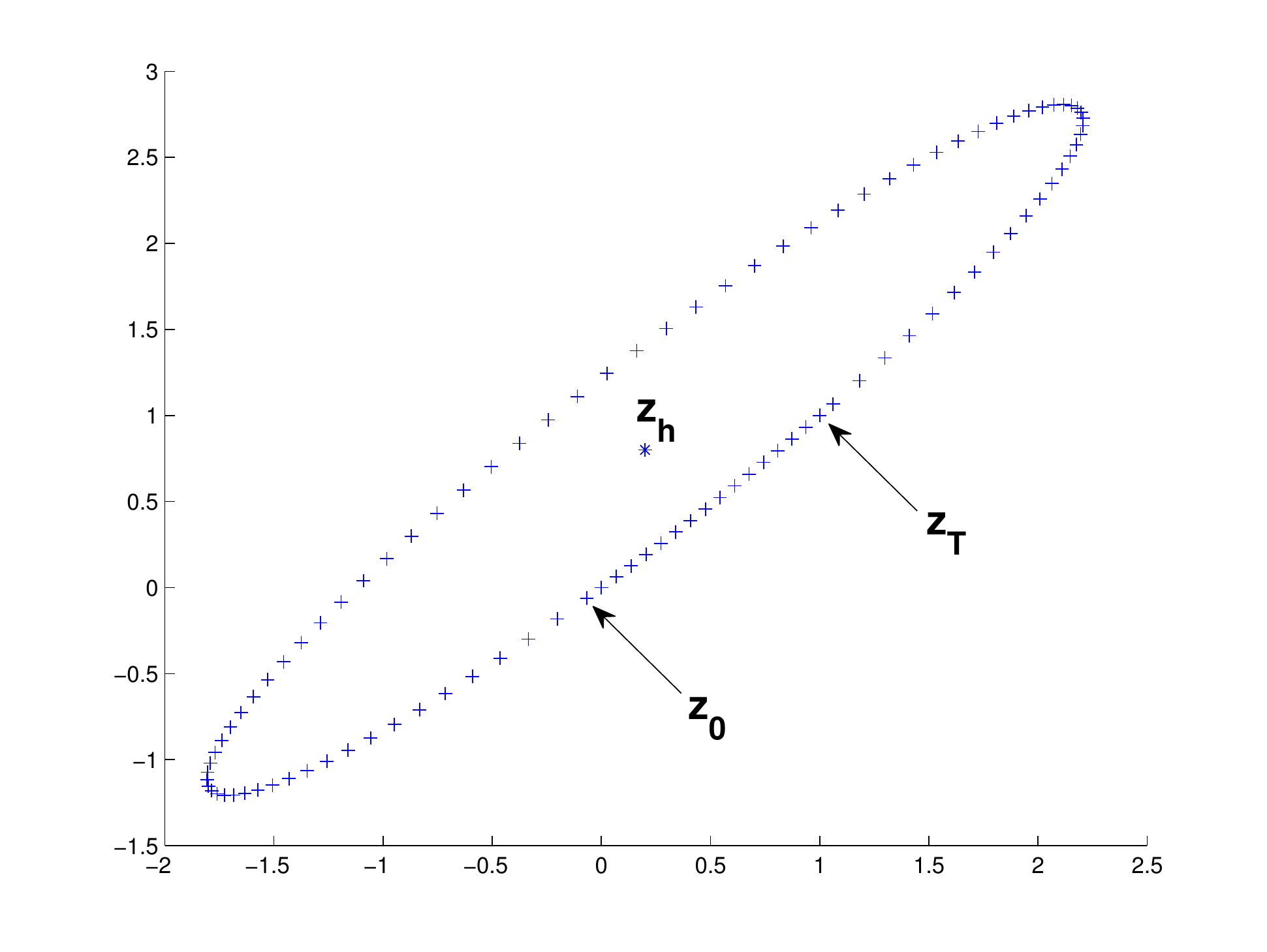}
        \caption{$T$ is larger than the ellipse period.}
        \label{fig:u0positive-ellipse}
    \end{subfigure}
    \caption{Traffic hole ($u_0>0$). Circles are iso-traffic levels.}
    \label{fig:u0positive}
\end{figure*}

We assume that there is a hot spot or a traffic hole located in $z_h$ and that the traffic intensity is of the form $u(t,z)=\frac{1}{2}u_0||z-z_h||^2+u_1=\frac{1}{2}u_0||z||^2-u_0z\cdot z_h+\frac{1}{2}u_0||z_h||^2 +u_1$. We can apply Corollary~\ref{cor:quad} with $e=-z_h$.
Figure~\ref{fig:hotspotu0negative} shows optimal trajectories when $z_h$ is a hot spot, i.e., for $u_0<0$, and different values of $K$. The starting point is $z_0$ and the destination is $z_T$. When $K$ increases, the velocity cost increases and the trajectories tend to the straight line between $z_0$ and $z_T$, which minimizes the speed. When $K$ is small, the UAV can go fast to $z_h$, reduces its speed in the vicinity of the hot spot and then goes fast to the destination
(in order to decrease
the cost function 
\eqref {eq:cost}
by increasing its traffic contribution).
Figure~\ref{fig:u0positive} shows optimal trajectories when $z_h$ is a traffic hole, i.e., for $u_0>0$. In Figure~\ref{fig:u0positive2}, $T$ is smaller than the period of the ellipse, i.e., $\frac{2\pi}{\omega}>T$. 
When $K$ decreases, the UAV can spend more time in the areas of higher traffic intensity. 
In Figure~\ref{fig:u0positive-ellipse}, $T$ is larger than the period. In this case, the trajectory follows one period of the ellipse whose equation is given by (\ref{eq:ellipse}) plus a part of the same ellipse from $z_0$ to $z_T$. 

\subsection{Multi-Phase Trajectory Characterization} \label{sec:interface}

We now consider a traffic intensity (or potential) consisting in two quadratic functions separated by an interface $\mathcal{I}$ of equal potentials delimiting two regions $1$ and $2$. The interface is defined by an equation $f(z)=C$, where $C$ is a constant and $f$ is a differentiable function. 
We assume that the optimal trajectory crosses the interface
only once,
at position $\xi$ and time $\tau$.
The impulsion $p^*$ is defined 
everywhere on the optimal trajectory between $(t_0,z_0)$ and $(T,z_T)$.
The following notations will be used  in the sequel: 
\begin{equation}
\label{eq:2-phase-notations}
p^- = 
p^*(\tau^-)=
\lim\limits_{\substack{s\to \tau\\s<\tau}}\ p^*(s)
\enspace,
\
p^+ = 
p^*(\tau^+)=
\lim\limits_{\substack{s\to \tau\\s>\tau}}\ p^*(s)
\enspace,
\ 
H^- = H_1(\xi,p^*(\tau^-))
\enspace,
\
H^{+} = H_2(\xi,p^*(\tau^+)) .
\end{equation}
\begin{theorem} \label{th:multiphase}
The location and time $(\xi, \tau)$ of interface crossing are characterized by the following equations
\begin{eqnarray}
p^- - p^+
-\mu\ \nabla_z f(\xi)&=&0 
\label{impulsion-ms:eq}
\\
H^+ - H^- &=& 0
\label{eq:Hconserved} \\
f(\xi)&=&C \label{eq:interface}
\end{eqnarray}
for some Lagrange multiplier $\mu\in\mathbb{R}$ 
\end{theorem}
\begin{proof}
See Appendix~\ref{app:multiphase}.
\end{proof}

Equation \eqref{eq:Hconserved} expresses the fact the energy is conserved when crossing the interface. One can show that actually the energy is conserved along the whole trajectory. Equation \eqref{impulsion-ms:eq} is related to the conservation of the tangential component of the impulsion at the interface. Equation \eqref{eq:interface} is the interface equation at $\xi$. 
One can show that under the assumption of equal potential on the interface, the kinetic energy, the impulsion, and the velocity vector are conserved across the interface.


\subsection{Multi-phase optimal trajectory uniqueness
and value function convexity 
issues}

Looking for uniqueness/non-uniqueness of optimal multi-phase trajectories
leads naturally to study the convexity of
the multi-phase total cost.
As stated in Appendix 
\ref{app:multiphase}
this cost is additive
and satisfies the dynamic programming principle which reads
\begin{equation}
\bar{S}(\Theta = \tau,\ \Xi = \xi\ |\ t_0,z_0,T,z_T)
=S_1(t_0,z_0,\tau,\xi)+S_2(\tau,\xi,T,z_T)
\label{eq:sumcost} 
\end{equation}
where $S_1$ and $S_2$ are themselves minimal. 
Hamilton-Jacobi equations
applied to each cost component 
allows to study
their
first- an second- order differential properties
such as convexity.
For instance in the chosen quadratic model,
each single-phase (minimal) total cost
$S_1(t_0, z_0, \tau,\xi)$, 
$S_2(\tau,\xi, T, z_T)$
given by Corollary \ref{cor:quad}
and \eqref{eq:vu0neg}
is convex (since it is a positive quadratic form with respect to the spatial coordinates)
with respect to the spatial position $\xi$,
but not necessarily so with respect to the interface crossing time $\tau$.
%
We consider in what follows a general
form of the single-phase value function
between $(t_1,x_1)$ and $(t_2,x_2)$ that is denoted by
\[
S(t_1, x_1,t_2, x_2).
\]
We also note  
by
$\w$\  
the pulsation,
$\phase = \w\ (\ttwo - \tone)$
is
the {\em temporal} phase
and 
$\pone,\ \ptwo$
are the
the initial and final impulsions
as derived from 
formula \eqref{eq:speedhotspot}.

\begin{figure}[t]
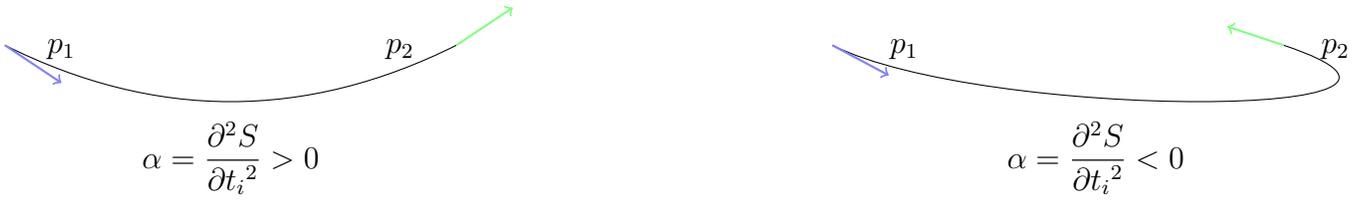
 
\tikzpicture
\centering
\coordinate (A) at (0,0) ;
\coordinate (B) at (6,0) ;
\draw
        (A) .. controls (2,-1) and (4,-1) .. (B) 
	node[very near start, above] {$\pone$} 
	
	node[very near end, above] {$\ptwo$} ;	
	
\draw [color=blue!50, thick] [->] (A) -- +(+0.75, -0.5); 
	
\draw[color=green!50, thick]  [->]  (B) --+ (0.75, 0.5) ; 

\draw  (A) + (3, -1.5) node 
{$
\alpha =\ddpart{S}{t_i} > 0$  } ;
\coordinate (C) at (0,3) ; 
\coordinate (D) at (6,3) ; 
\coordinate (C) at (11,0) ; 
\coordinate (D) at (17,0) ; 

\draw
        (C) .. controls +(2,-1) and +(3,-1) .. (D) 
	
	node[very near start, 
	above] {$\pone$} 
	node[very near end, 
	above] {$\ptwo$} ;	
	
\draw [color=blue!50, thick]  [->] (C) -- +(+0.75, -0.4); 
\draw[color=green!50, thick]  [->] (D) -- +(-0.75, 0.25) ; 

\draw (C) + (3.5, -1.5) node 
{$
\alpha =\ddpart{S}{t_i} < 0$} ;
\endtikzpicture
\caption{
{Physical} interpretation of local convexity vs. non-convexity
 of the single-phase cost:
 \\
 the initial and final impulsions
 $(\pone, \ptwo)$ form an {\em obtuse} angle 
 $\Ra$ the value function is {\em non-convex}. 
 }
\label{fig:p1-p2-representation}
\end{figure}

\begin{theorem} \label{th:single-phase-hessian}
i) The Hessian of the single-phase cost 
\wrt 
any joint variable
$\psi_i = $
${(T_i, X_i)}_{i=1,2}$ 
\ 
is
\begin{align}
\label{full-single-phase-hessian:eq}
{\cal H}(\psi_i)
=
\nabla^2_{T_i, X_i}\ S(t_1, x_1,t_2, x_2)
&=
\Hessian{\alpha}{K g}{\Pi_i}
\hspace*{1mm} 
\mbox{ with }
\hspace*{1mm} 
\left\{
\begin{array}{l}
\alpha =  
\dfrac{\partial^2 S}{\partial T_i^2}
= \dfrac{\w\ \pone \cdot \ptwo}{K\ \sinh \phase}
\\[4mm]
g = \w\ \coth \phase\ > 0 
\\[4mm]
\Pi_i 
= - \dfrac{\w\ p_{3-i}}{\sinh \phase} \in \Real^2
\end{array}
\right.
\end{align}
where  $p_i$ 
is the impulsion at 
time $T_i$, \ie at extremity $X_i$
\ ($i = 1,2$ and $j = 3-i$).
\\
ii)
At most one eigenvalue 
of ${\cal H}(\psi_i)
$
can be  negative 
and $\alpha < 0$ 
is a 
sufficient condition for this to hold,
namely implying local non-convexity of the 
(single-phase)
value function.
\\
iii) 
The 
double-phase
Hessian enjoys a similar structure
as 
\eqref{full-single-phase-hessian:eq}
and writes
with respect to the variable 
$ 
\Psi =
(\Theta, \Xi)$
\begin{subequations}
\begin{align}
\label{eq:Hessian-second-a}
\bar{{\cal H}}(\Theta =\tau, \Xi = \xi) 
=
& %
\ 
{\cal H}_1(\chi_2) + {\cal H}_2(\chi_1) 
=
\Hessian{\alpha}{K h}{\Pi}
\\[2mm]
\label{eq:Hessian-second-b}
\mbox{with }
& %
\left\{
\begin{array}{l}
\alpha = 
\ddpart{S_1}{\Ttwo}(t_0,z_0,\tau,\xi)
+ 
\ddpart{S_2}{\Tone}(\tau,\xi,T,z_T)
\ \
(\mbox{cf. } \eqref{full-single-phase-hessian:eq})
\\[4mm]
h = \w_1\ \coth \phase_1 + \w_2\ \coth \phase_2
\\[4mm]
\Pi = 
- \left[ 
\dfrac{\w_1\ p(t_0)}{\sinh\phase_1}
+ \dfrac{\w_2\ p(T)}{\sinh\phase_2}
\right]
\\[4mm]
\mbox{and } 
\phase_1 = \w_1\ (\tau - t_0),\
\phase_2 = \w_2\ (T- \tau )
\end{array}
\right.
\end{align} 
\end{subequations}
and may thus be non-convex
as well.
\end{theorem}
\begin{proof}
See Appendix \ref{app:Hessians}
subsections
\ref{structure-hessian:app}
and
\ref{diagonal-hessian:app}. 
\end{proof}
In the single-phase case,
the sufficient condition 
{\em ii)} 
\ie
$\alpha < 0$,
supports first a {\em physical}
interpretation
(see Fig. \ref{fig:p1-p2-representation})
and also a {\em geometrical} one
(see next 
Theorem \ref{th:two-phase-convexity} 
and Fig. \ref{fig:u-v-representation}).
\begin{theorem}
Let $u=\frac{x_1+x_2}{2}$ and $v=\frac{x_2-x_1}{2}$. (Note that $\|u\|$ is the distance from the hotspot to the average of $x_1$ and $x_2$, and the $\|v\|$ is the half-distance between $x_1$ and $x_2$.)  
\label{th:two-phase-convexity}
Then the following sufficient condition holds
\[
\dfrac{\parallel v \parallel}{\parallel u \parallel} < 
\dfrac{\cosh \phase - 1}{1 + \cosh \phase}
= \tanh^2 \dfrac{\phase}{2} 
\Rightarrow
\alpha =\ddpart{S}{t_i} < 0 .
\]
\end{theorem}

\label{app:geometrical-non-convexity-condition}
{\em Proof}:
We combine the 
formulas 
\eqref{pone:eq} and \eqref{ptwo:eq}
in 
\eqref{d2S-dphase2:eq}
to obtain 
\begin{equation}
\alpha = \dfrac{w}{K}\ 
\dfrac{\pone \cdot \ptwo}{\sinh \phase}
=
\dfrac{K\,\w^3} 
{\sinh^3 \phase}
 \left[ (\parallel \xone\parallel^2 + \parallel \xtwo\parallel^2)\ \cosh \phase - \xone\ .\ \xtwo\ (1 + \cosh^2 \phase) \right]
 \label{alpha=d2Sdt2:eq}
\end{equation}
which can be rewritten as
\[
\alpha = \dfrac{K\,\w^3} 
{\sinh^3 \phase}
\left[ - \parallel u \parallel^2\ ( 1 - \cosh \phase)^2 + 
\parallel v \parallel^2\ (1 + \cosh \phase)^2)
\right].
\hfill
\qed
\nonumber
\]
This non-convexity condition for the single-phase value function 
interprets as:
long phase 
($\tanh\frac{\phase}{2} > 
\frac{1}{2}$)
and long distances between the hotspot $z_h$
to both
initial and final positions $\xone$, $\xtwo$,
relatively to their mutual distance
($\frac{\parallel u \parallel}
{\parallel v \parallel} \ggg 1$).
\begin{figure}[t] 
~\\[-10mm]
\centering
\begin{tikzpicture}[scale = 0.6, information text/.style={rounded corners, inner sep=1ex}]

\coordinate [label=above: 
\textcolor{blue}{$z_h$}] (htspt) at (-4,2) ; 

\coordinate [label=left: \textcolor{blue}{$\xone$}] (x1) at (-1,-2) ;

\coordinate [label=right: \textcolor{blue}{$\xtwo$}] (x2) at (+1,-1) ;

\draw[->] (x1) -- (x2) ; 

\coordinate 
(m) at ($(x1) !.5! (x2)$) ; 

\coordinate [label=above: 
\textcolor{purple}{$v$}] () at ($(m) !.5! (x2)$) ;   

\draw [->]  (htspt) to [edge label = \textcolor{purple}{$u$} ] (m) ;

\draw[xshift=2.0cm, yshift=+2cm]
node[right, text width=6cm, information text] {
\begin{align*}
u &=  \dfrac{\xone + \xtwo}{2} 
\\
v &= \dfrac{\xone - \xtwo}{2} 
\end{align*}
} ; 
\end{tikzpicture}
\caption{ {Geometrical}
interpretation of local  non-convexity:
long time
and long distance to the hotspot.
}
\label{fig:u-v-representation}
\end{figure}
Several preliminary simulations
show that when this case 
happens for two phases, the total value function is indeed non-convex
and that 
several (local) optimal solutions may indeed exist.
\section{Algorithms} \label{sec:algos}
Previous propositions
allow us to propose several numerical algorithms
seeking optimal trajectories.
In this section, we present two algorithms aiming this goal:
a gradient descent algorithm 
and
a bisection search method based
on the linear control property
\eqref{eq:speedhotspot}.

\subsection{A Gradient Descent Algorithm  \GRADALGO}
In this section, we propose \textsc{Grad-Algo} (Algorithm~\ref{alg:ms-grad-algo}) which is an alternated optimization-based algorithm relying on the following procedures. 

\paragraph{Procedure for seeking an optimal $\xi$ given a fixed  $\tau$}
from Hamilton-Jacobi equations
(Lemma \ref{lemma:hamiltonjacobi}
and Appendix \ref{app:hamiltonjacobi})
the gradient of the total cost function with respect to $\xi$ is
$
p^*(\tau^-) - p^*(\tau+)$.
Equation~\eqref{eq:impulsionquad} 
says that 
in the Newtonian framework the impulsion is proportional to the control variable $a$.
Equation~
\eqref{eq:speedhotspot}
says that in the quadratic model the velocity vector is, at any time a linear combination of {\em centered} initial and final positions.
Then
the searched gradient appears to be an {\em affine} function of $\xi$ which reads
\begin{align} 
p^*(\tau^-) - p^*(\tau+) 
&= 
\nabla_{\Xtwo}
S_1(t_0,z_0,\tau,\xi)
+
\nabla_{\Xone}
S_2(\tau,\xi,T,z_T) \
\nonumber
\\
& = 
K h\ (\xi - B).
\label{grad-impulse-difference:eq} 
\end{align} 
Scalar $h$ and point $B$ (where the spatial gradient cancels 
 at fixed $\tau$
\ie\ $p^*(\tau^-) = p^*(\tau+)$)
verify
\begin{eqnarray}
h &=& \omega_1\, \coth(\omega_1(\tau -t_0)) + \omega_2\, \coth(\omega_2(T - \tau)),
\label{h-coefficient:eq} 
\\
B &=&
\dfrac{1}{h}\big[\
\omega_1\, \HSPONE\, \coth(\omega_1 (\tau - t_0)) 
+ \omega_2\, \HSPTWO\, \coth(\omega_2(T - \tau)) 
\notag \\
&&+\dfrac{\omega_1\, (z_0 - \HSPONE )}{\sinh(\omega_1 (\tau - t_0))} 
   + \dfrac{\omega_2\,  (z_T - \HSPTWO)}
      {\sinh (\omega_2 (T - \tau))} 
      \label{location-B-ms:eq}
\ \big].
\end{eqnarray}
The equation involving the Lagrange multiplier 
(\ref{impulsion-ms:eq}) now reads
\begin{eqnarray}
K h\ (\xi - B) - \mu\ \nabla_{\xi} f(\xi) = 0
\label{algo-proj-B:eq}
\end{eqnarray}
and shows that the optimal location $\xi^*$ 
is the {\em orthogonal projection} of 
$B$
on the interface 
${\cal I}$.

\paragraph{Procedure for seeking an optimal $\tau$ given a fixed  $\xi$}
we use the result of Theorem~\ref{th:multiphase}. 
As also shown by Hamilton-Jacobi equations
(Lemma \ref{lemma:hamiltonjacobi}
and Appendix \ref{app:hamiltonjacobi}),
the gradient of $S$ with respect to $\tau$ is given by $H_2(\xi,p^*(\tau^+))-H_1(\xi,p^*(\tau^-))$. The Hamiltonians 
are easily computed in each phase by applying
the classical Newton formula at given location
$(z,t)$
\[
H(p, z, t)  = 
\dfrac{{\norm{p}}^2}{2\ K} + V(z).
\]
We then update $\tau$ by using a simple gradient descent scheme
(see Algorithm \ref{alg:ms-grad-algo}).

\def\abs#1{| {#1} |}
\def\cP{cP} 
\def\cH{cH} 
\paragraph{Stop criteria} 
the search procedure is performed until the impulsion and the Hamiltonian have converged with a given accuracy. For this purpose  we consider the function 
$
g:(\Real^n \setminus \{0\}) \times (\Real^n \setminus \{0\}) 
\mapsto \Real
$
\\
defined as follows
$g(x, y) = \dfrac
{\norm{x-y}}      {\inf(\norm{x},\norm{y})}
\ \ \ 
\mbox{ which measures a relative 
``discrepancy"
between 
$x$ and $y$}.
$
Now let us give 
$\epsilon_p$ and $\epsilon_H$ 
(typically of the order of $2.10^{-4}$)
and define the stop criteria as follows
\begin{equation}
 \texttt{STOP}  =
g(p^+, p^- ) ~/~ \epsilon_p
< 1 
\ \ \
\&\& 
\ \
g(H^{+}, H^{-}) ~/~ \epsilon_H 
< 1.
\label{eq:stop-criteria}
\end{equation}
\def\myproj{\mbox{Proj}_{{\cal I}}}
 \begin{algorithm}[t] 
\caption{\textsc{Grad-Algo} An 
uncoupled 
projected gradient descent algorithm}
\label{alg:ms-grad-algo}
\begin{algorithmic}[1]
\State {\bf Input:}
precisions $\epsilon_p$, $\epsilon_H$ 
($\epsilon_p = \epsilon_H$ 
= $2.10^{-4}$), 
$M_{\tau}=$
number of gradient descent iterations
on $\tau$
\vspace*{1mm}
\State {\bf Init:}
starting position
$\xi \gets$
     $\xi_0$ 
     e.g., 
     = $\frac{z_0 + z}{2}$
     ,
 $\tau \gets$ 
 $\tau_0$ 
 e.g., 
 = $\frac{t_0 + T}{2}$
\State {\bf Output:} 
{
$(\tau, \xi)$}
\def\BIPHASE{\texttt{TWO\_PHASE}}
\Procedure{\BIPHASE}{{
$\tau, \xi$}}
    \State {\bf computes:} 
the current bi-phase trajectory given by
$(t_0, z_0) \ra (\tau, \xi) \ra (T, z_T)$
\State \Return
$(B, H^+, H^-, 
\bar{\cal H}
)$ 
= 
(the 
B-point
\eqref{location-B-ms:eq},
the two phase Hamiltonians
and 
{
the total Hessian} 
at 
current location 
$(\tau, \xi)$
\
(\ref{eq:sumcost},
\ref{eq:Hessian-second-a},
\ref{eq:Hessian-second-b}))
\EndProcedure
\State $N_{step} = 0$
\State
$(B, H^+, H^-, 
{
\bar{\cal H}})
=
\BIPHASE(
{
\tau, \xi})$
\Repeat
\State
$\xi \gets \myproj(B)$ 
\eqref{algo-proj-B:eq}
\vspace*{1mm}
 \For{$m=1,\dots,M_{\tau}$}
 \vspace*{1mm}
 \State
 $\tau \gets \tau - \dfrac{(H^+ - H^-)}{\alpha}$
 \ 
 with 
 $\alpha 
= {\bar{\cal H}}_{11}
= \dfrac{\partial^2 {\bar S}}{\partial \Theta^2}(\tau,\xi)
$
\ 
 (Newton descent 
 \wrt $\tau$)
 \vspace*{1mm}
 \State
 $(B, H^+, H^-, 
 {
 \bar{\cal H}}
 ) = 
 \BIPHASE(
 {
 \tau, \xi})$
 \EndFor
    \State $N_{step} = N_{step} + 1 $
    \Until \texttt{STOP} \eqref{eq:stop-criteria}
\end{algorithmic}
\end{algorithm}
{\bf Complexity:}
Algorithm
\ref{alg:ms-grad-algo}:
if 
$\alpha^* = \ddpart{\bar{S}}{\Theta}(\tau^*,\xi*)
\neq 0
$
then the convergence of Newton descent
on $\tau$ is quadratic (see~\cite{bierlaire.15.book} for instance).
If $\alpha^* = 0$ 
say,
$
H_2(\xi,p^*(\tau^+))-H_1(\xi,p^*(\tau^-))
\sim
(\tau-\tau^*)^\nu$
with $\nu \geq 2$,
the convergence is linear (see~\cite{bierlaire.15.book} for instance).
For safety we take
$M_\tau = 10$
\ie
$\nu \approx 2$.
The global complexity  is
$O(N_{step}\ M_{\tau})$
and usually $N_{step} = 2$.

\subsection{The B-curve algorithm
\BALGO}
\label{B-curve:sec}
The $B$-curve algorithm 
aims to
overcome the
non-convexity issues developed 
previously. 
It proceeds as follows.
Consider the $B$ point defined in formula
\eqref{location-B-ms:eq}.
We notice that for each 
$\tau \in [t_0, T]$ 
this point is defined univocally knowing all parameters 
($t_0,\ T,\ z_0,\ z_T,\ z_h$)
and all (quadratic) traffic profiles
so that we can see $B$ in \eqref{location-B-ms:eq} as a function of $\tau$, 
denoted
$B(\tau)$.
It can be shown that
this function 
$\tau \mapsto B(\tau)$ 
is continuous 
in the interval $[t_0, T]$ 
and that $\lim_{\tau \ra t_0^+}{B(\tau)} = z_0$,
$\lim_{\tau \ra T^-} B(\tau) = z_T$.
Assume now that the two optimal 
sub-trajectories
are such that crossing time $\tau \in [t_0, T]$
and interface position $\xi$ verify
$\xi = B(\tau)$.
Then by formula 
\eqref{grad-impulse-difference:eq}
the spatial gradient of $S$ at point $\xi$ is
\begin{equation}
    p^*(\tau^-) - p^*(\tau+) = 0 
    \label{eq:delta-p-B-algo}
\end{equation}
This implies that the  {\em kinetic} components
of both  Hamiltonians are equal at the interface.
Now since $\xi$ belongs to the interface,
both phase {\em potentials} (traffics) 
are 
equal
by definition,
so that the {\em total} Hamiltonian is also conserved at the interface:
this is 
the optimality condition required
with respect to time $\tau$.
To summarize,  
in these conditions,
both the Hamiltonian and the total impulsion are conserved
at interface $\mathcal{I}$, i.e., local optimality conditions hold for the total value function.
It is also worth noticing that the related Lagrange multiplier appearing in 
\eqref{algo-proj-B:eq} 
now simply vanishes: $\mu = 0$.
The proposed  B-algorithm
consists then in seeking the intersection
of the B-curve 
${\cal B} = \{ B(\tau), \tau\in[t_0,T]\}$
with the interface ${\cal I}$
(Algorithm \ref{alg:bcurve-dichotomy}). 
For this, 
we first select precision $\epsilon_B$,  
then proceed by bisection
and check the stopping criterion 
\def\STOPB{\texttt{STOP}\_\texttt{B}}
\begin{equation}
\STOPB
= \dfrac{\abs{t_2 - t_1}}{\abs{b - a}} 
\ \dfrac{1}{\epsilon_B}
\ \ \ 
< 1.
\label{eq:stop-B}
\end{equation}
{\bf Complexity:} Algorithm
\ref{alg:bcurve-dichotomy}:
if $\epsilon_B = \dfrac{1}{2^m}$
then the bisection algorithm converges in $m$ iterations. 
Its complexity is thus
$O(\log   \dfrac{1}{\epsilon_B})$.
In our experiments
we chose  
$\epsilon_B =2.10^{-4}$ 
$\approx \dfrac{1}{2^{12}}$
and 
12 iterations 
are 
indeed 
sufficient to provide a trajectory with optimality conditions holding
at this (relative) precision. 
\def\MYBOOL{\texttt{BOOL2}}
\def\MYBOOL{\texttt{IN\_ZONE2}}
\begin{algorithm}[t]
\caption{the B-curve bisection algorithm}
\label{alg:bcurve-dichotomy}
\begin{algorithmic}
\State {\bf Input:} initial/final time 
search interval $[a,b] \in [t_0, T]$,
precision $\epsilon_B$
fixed by the user
($\epsilon_B =2.10^{-4}$) 
\State {\bf Output:} $(\xi, \tau)$
\State {\bf Init: } $t_1\gets a \enspace, t_2\gets b$,
the algorithm only starts if \MYBOOL($B(a)$) $\neq$\ \MYBOOL($B(b)$)
\Procedure{\MYBOOL}{$z$}
\State
$u_i(z)$ = (time-stationary) traffic $u(z)$ 
from hotspot $z_{hi}$ at point $z$
\ ($i = 1,2$) 
\State \Return $(u_2(z) > u_1(z))$
\EndProcedure
\vspace*{1mm}
\Repeat
    \State $x_1 \gets B(t_1)$
    \enspace, 
    $x_2 \gets B(t_2)$ 
    \vspace*{1mm}
    \State
     $\tau \gets \dfrac{t_1 + t_2}{2}
     \enspace,\ \xi \gets B(\tau)$
     \vspace*{1mm}
     \If{(\MYBOOL($\xi$) == \MYBOOL($x_2$))}
		\State $t_2 \gets \tau$ 
	  \Else \ \ $t_1 \gets \tau$
	\EndIf
\State {\bf compute stop criterion}
\STOPB
\ \ 
\eqref{eq:stop-B}
 \Until \STOPB
\vspace*{1mm}
   \State $\xi \gets B(\tau)$
\end{algorithmic}
\end{algorithm}

\section{Numerical Experiments} \label{sec:numericalexp}

\subsection{From Measurements to Quadratic Profile}

\begin{figure*}[t!]
    \centering
    \begin{subfigure}[t]{0.5\textwidth}
        \centering
        \includegraphics[width=\linewidth]{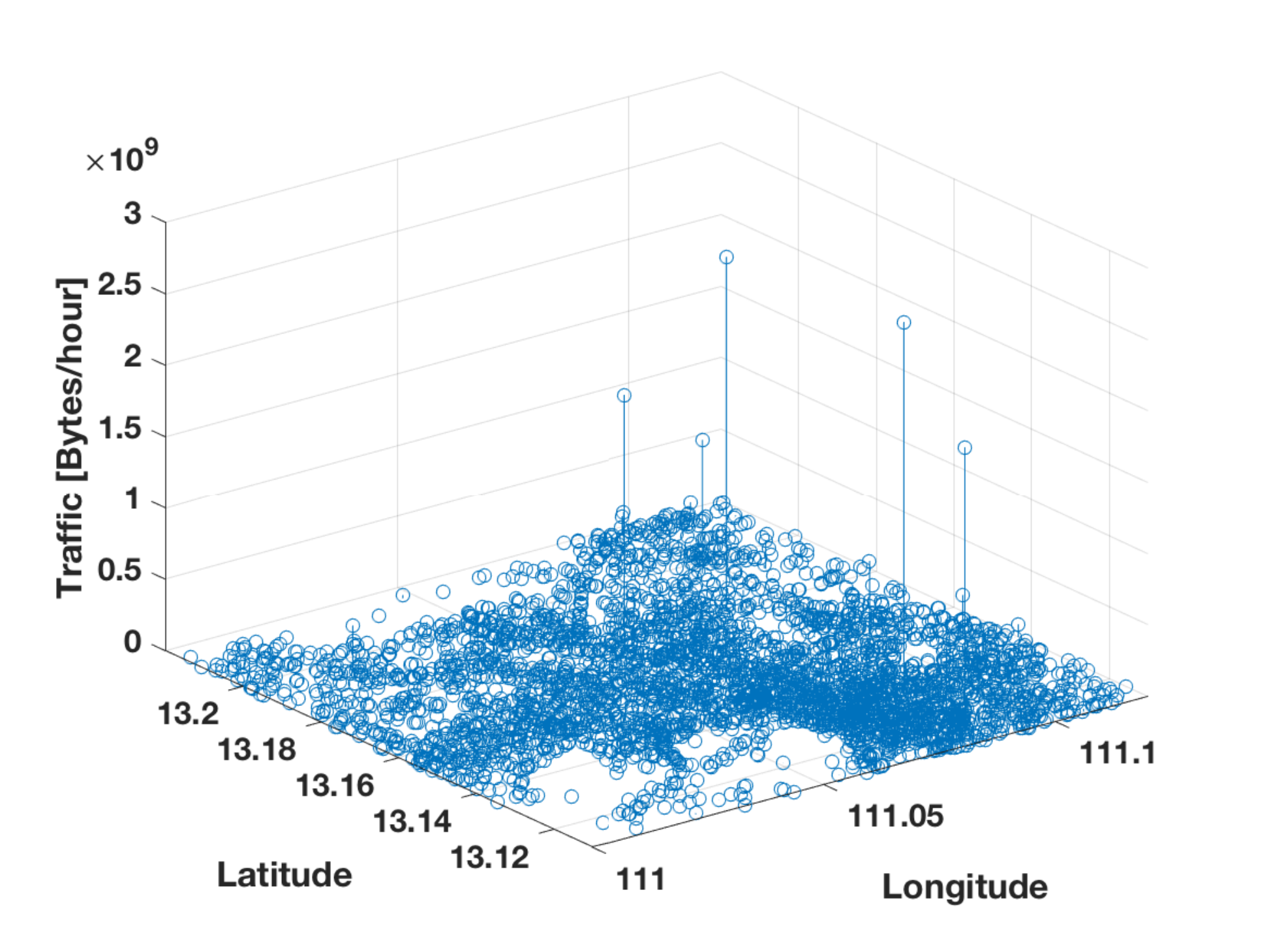}
        \caption{Raw traffic data on the 22th Aug. 2012 6pm \cite{Modeling15-Chen}. }. 
        \label{fig:rawdata}
    \end{subfigure}%
    ~ 
    \begin{subfigure}[t]{0.5\textwidth}
        \centering
        \includegraphics[width=\linewidth]{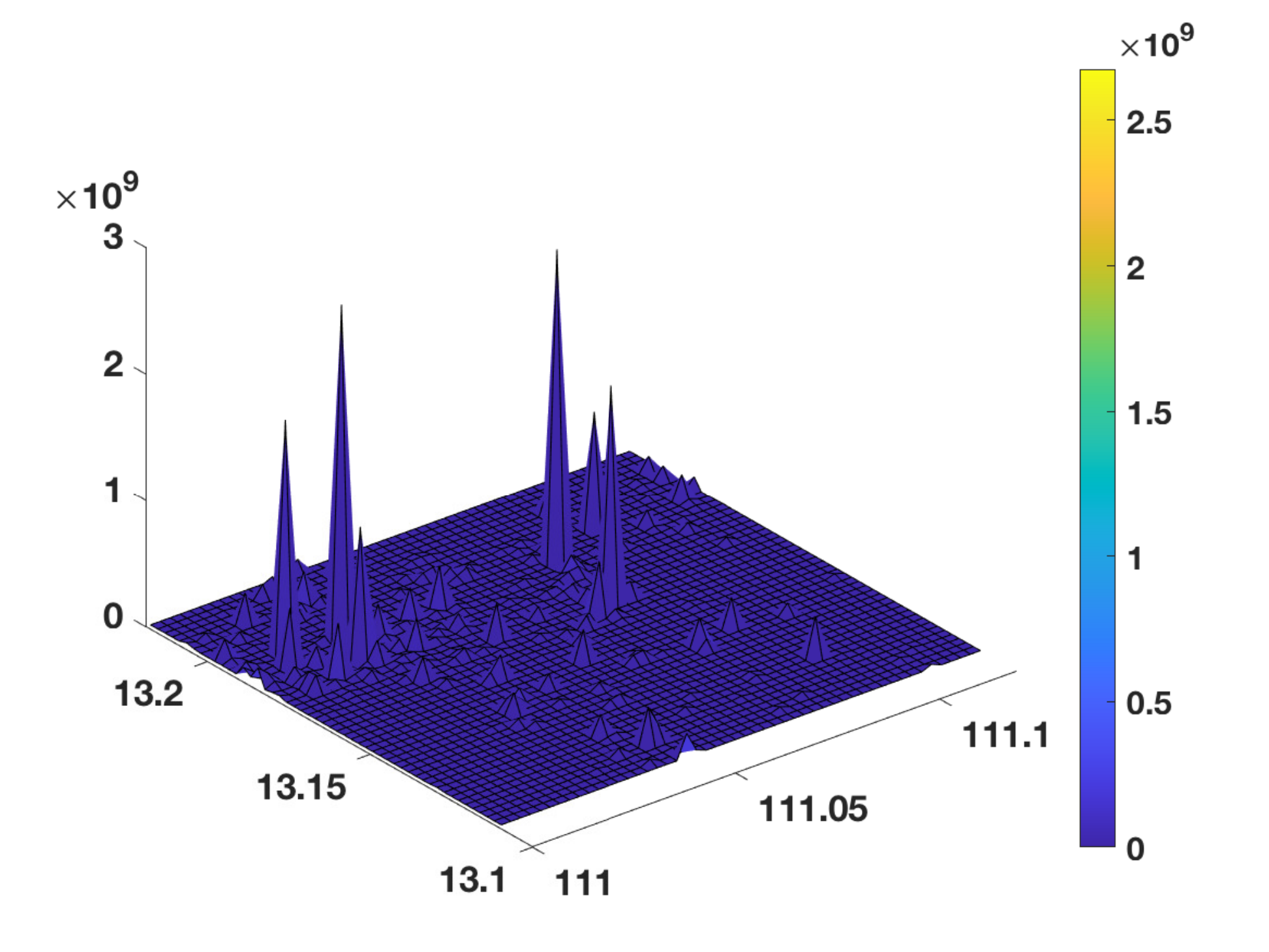}
        \caption{First smoothing: data aggregation at drone coverage level.}
        \label{fig:firstsmoothing}
    \end{subfigure}
    \caption{Data preprocessing: raw data and first smoothing.}
\end{figure*}

\begin{figure*}[t!]
    \centering
    \begin{subfigure}[t]{0.5\textwidth}
        \centering
        \includegraphics[width=\linewidth]{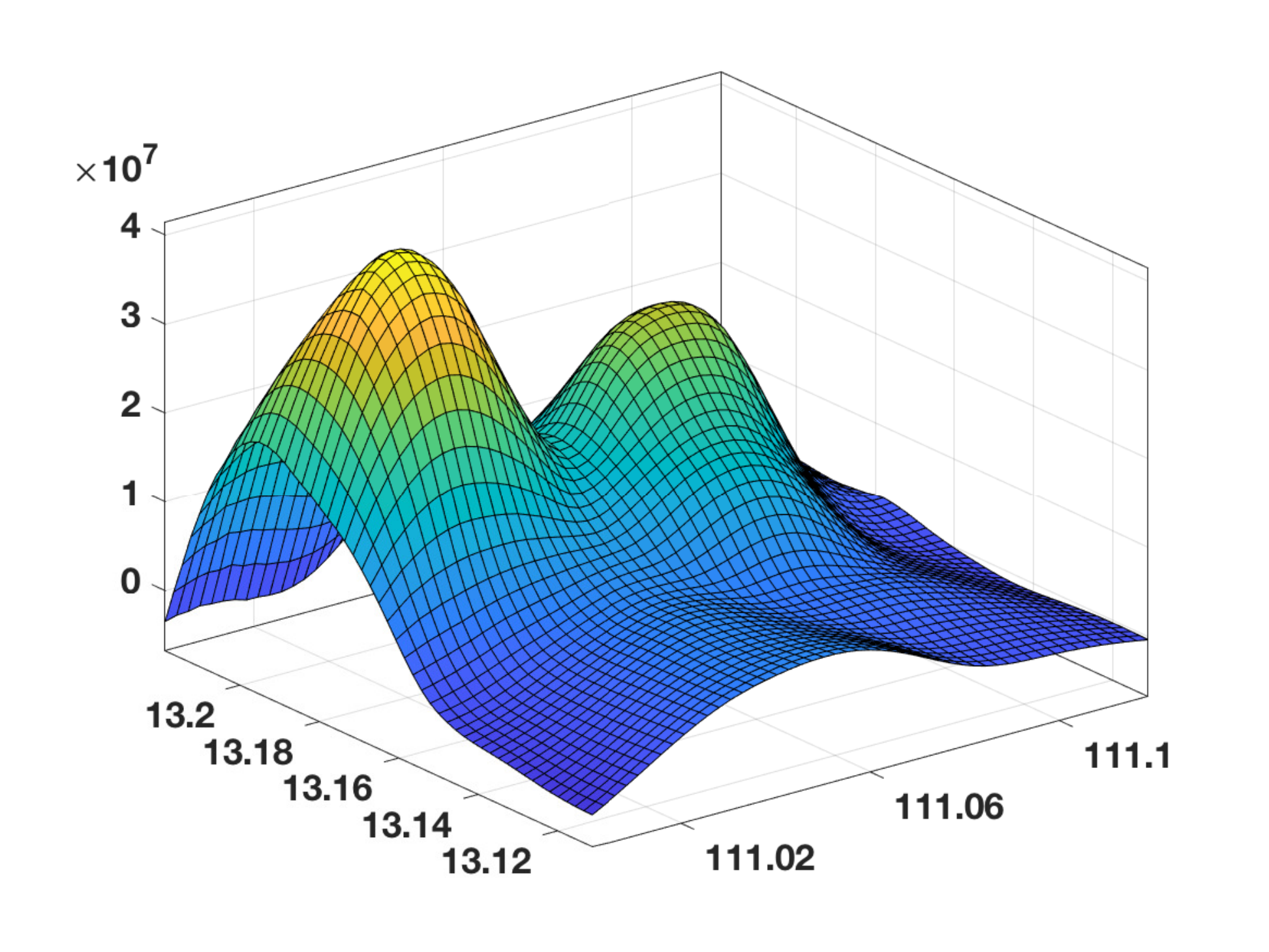}
        \caption{Second smoothing: LOWESS with $\alpha=0.25$.}
        \label{fig:secondsmoothing}
    \end{subfigure}%
    ~ 
    \begin{subfigure}[t]{0.5\textwidth}
        \centering
        \includegraphics[width=\linewidth]{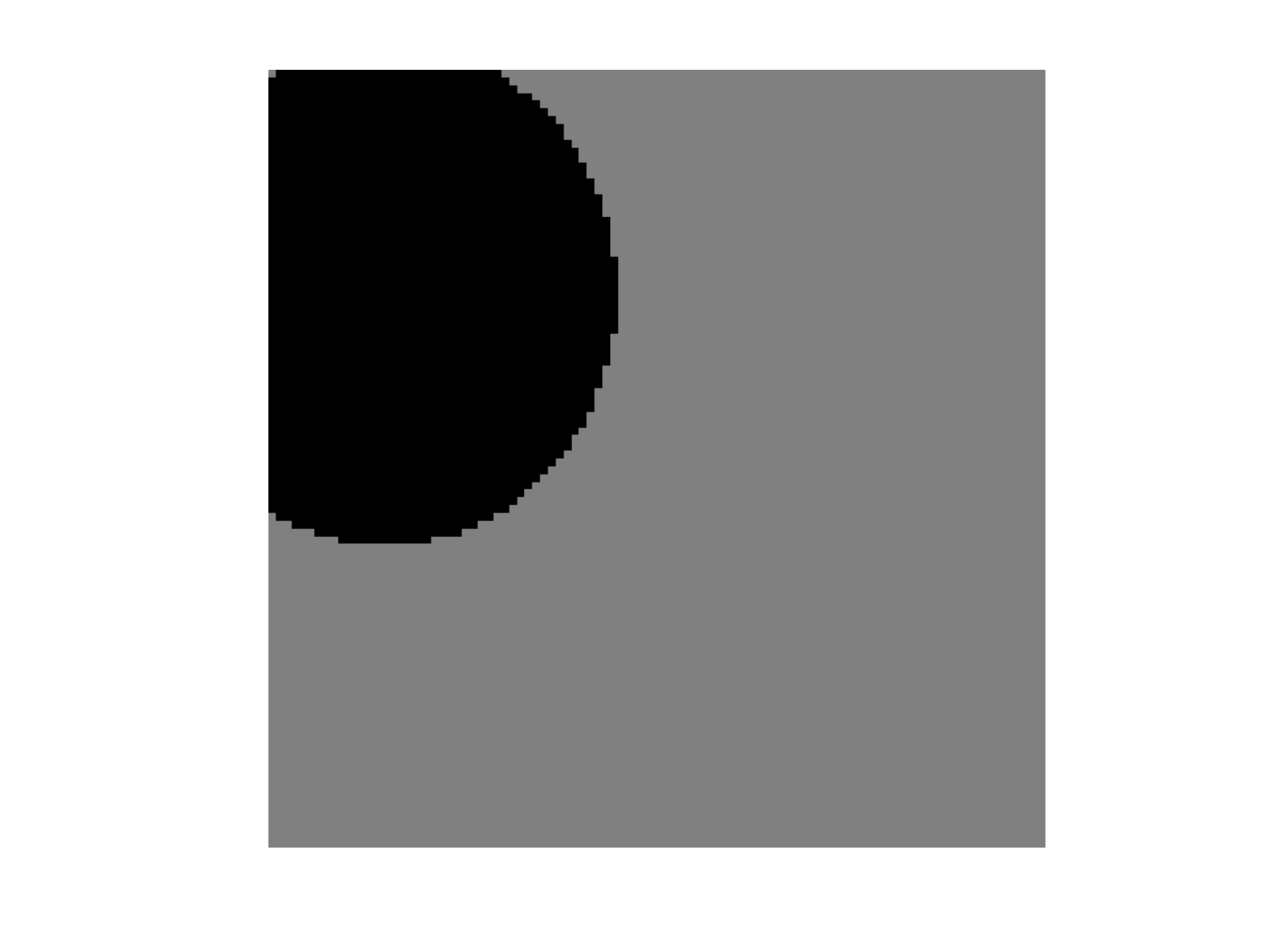}
        \caption{K-means partition.}
        \label{fig:kmeanspartition}
    \end{subfigure}
    \caption{Data preprocessing: second smoothing and partitions in 2 phases.}
\end{figure*}

In this section, we explain how from measured or estimated traffic load, we can derive a quadratic model that will allow us to apply our framework. To illustrate the procedure, we extract data from the open data set presented in~\cite{Modeling15-Chen}. Traffic data (in number of bytes) has been collected from an operational cellular network in a medium-size city in China and is recorded for every base station and every hour. For our experiment, we extract a rectangle region $[X_{\min},X_{\max}]\times [Y_{\min},Y_{\max}]$, where $X_{\min}=111$, $X_{\max}=111.12$, $Y_{\min}=13.12$, $Y_{\max}=13.22$ are the minimum and maximum longitude and latitude respectively (real figures have been anonymized). This corresponds approximately to a rectangle of $11$~km$\times 13$~km with $400$ base stations having an average cell range of $337$~m. The traffic of the 22th August 2012 between 5 and 6pm is illustrated in Figure~\ref{fig:rawdata}. We assume that this traffic is representative of the traffic intensity when the drone is launched. In order to fit this raw data to our model, we follow the pre-processing steps shown in Algorithm~\ref{alg:preproc}.
\begin{algorithm}[t]
\caption{Data preprocessing}
\label{alg:preproc}
\begin{algorithmic}[1]
    \State First smoothing: data aggregation
    \State Second smoothing: LOWESS
    \State Normalization
    \State K-means with quadratic models 
\end{algorithmic}
\end{algorithm}
The first smoothing consists in aggregating the traffic data on a grid of $50$ steps in both longitude and latitude directions. The resulting elementary regions should correspond approximately to the drone coverage. The result is shown in Figure~\ref{fig:firstsmoothing}. The data exhibits a very high variability with very high peaks around few locations. The second smoothing is a Locally Weighted Scatterplot Smoothing (LOWESS)~\cite{cleveland1979robust}. We use here the Matlab function \texttt{fit} with the option "Lowess". The choice of the smoothing parameter $\alpha$, i.e., the proportion of data points used for every local regression, has a decisive impact on the result. Increasing $\alpha$ has the effect of averaging out the different peaks. In our specific scenario, $\alpha=0.25$ yields  Figure~\ref{fig:secondsmoothing} with two local maxima. With $\alpha=0.5$, we obtain a single maximum. In step 3 of the pre-processing, the traffic is normalized between $0$ and $1$ with no influence on the optimal trajectory.  

The final pre-processing step is an adaptation of the classical K-means algorithm (see Algorithm~\ref{alg:kmeans}) to fit to quadratic models. Inputs are the data points obtained after the normalization, $K_c$, the number of clusters (or hot spots), $K_n$ the number of nearest neighbors and $M$ the number of iterations. Every cluster is associated to a quadratic function (in our case, we have $K_c=2$). Every data point $j$ is associated to a cluster and has a related label $L_j$ in $\{1,...,K\}$. For every data point, a list of nearest neighbors is built (step 4). An arbitrary initial labelization is chosen (step 6). The algorithm then proceeds by iterations (steps 7-16). At every iteration, if a point $j$ has some neighbor with a different label (step 9), a best new label is found for $j$ (step 10-13) in terms of quadratic error $e_k$. The error $e_k$ measures the difference between the data points and the $K_c$-quadratic model, which fits a quadratic function to every cluster assuming that $j$ has label $k$ (steps 18-31). In Algorithm~\ref{alg:kmeans}, $\texttt{KNN}(K_n,X,Y)$ is a procedure that finds the $K_n$ nearest neighbors of $(X,Y)$ with respect to the Euclidian distance. We use the Matlab implementation \texttt{knnsearch}.  $\texttt{NLLS}(\mathcal{L})$ is a non-linear least square method that fits data points in $\mathcal{L}$ to a quadratic function of the form $\frac{1}{2}u_{0}|z-z_{h}|^2+u_{1}$. The Matlab implementation is based on a trust-region approach of the Levenberg-Marquardt Algorithm.

\begin{algorithm}
\caption{K-means with quadratic models}\label{alg:kmeans}
\begin{algorithmic}[1]
\State {\bf Input:} $K_c$ (number of clusters), $J$ (number of measurement points) $(X_j,Y_j,Z_j)$, $j=1,\dots,J$ (coordinates and estimated traffic load $Z_j$ in $(X_j,Y_j)$), $K_n$ (number of nearest neighbors), $M$ (number of iterations)
\State {\bf Output:} $z_{hk}$, $u_{0k}$, $u_{1k}$ $k=1,\dots,K_c$ (hot spot characteristics), $e$ (quadratic error)
\For{$j=1,\dots,J$}
    \State $\mathcal{K}_j\gets\texttt{KNN}(K_n,X_j,Y_j)$
\EndFor
\State $L\gets$ Initial labelization 
\For{$m=1,\dots,M$}
    \For{$j=1,\dots,J$}
        \If{$\exists j'\in \mathcal{K}_j$ s.t. $L_{j'}\neq L_j$}
            \For{$k=1...K_c$}
                \State $e_k,z_{hl},u_{0l},u_{1l}$, $l=1,\dots,K_c$ 
                $\gets \texttt{FIT}((X,Y,Z),j,L,k,K_c)$
            \EndFor
            \State $L_j\gets \arg\min_k e_k$
        \EndIf
    \EndFor
\EndFor
\State \Return $z_{hk}$, $u_{0k}$, $u_{1k}$ $k=1,\dots,K_c$, global quadratic error
\Procedure{\texttt{FIT}}{$(X,Y,Z)$, $j$,$L$,$k$,$K_c$}
    \State $L_j\gets k$
    \For{$l=1$,...,$K_c$}
        \State $\mathcal{L}_l \gets \{(X_i,Y_i,Z_i)|L_i=l\}$
        \State $z_{hl},u_{0l},u_{1l}\gets \texttt{NLLS}(\mathcal{L}_l)$
    \EndFor
    \State $e_k\gets 0$
    \For{$j=1$,...,$J$}
        \State $\tilde{Z}_j \gets \max_{l=1,\dots,K_c}\frac{1}{2}u_{0l}|z-z_{hl}|^2+u_{1l}$
        \State $e_k\gets e_k+|\tilde{Z}_j-Z_j|^2$
    \EndFor
    \State $e_k\gets e_k/J$
    \State \Return $e_k,z_{hl},u_{0l},u_{1l}$, $l=1,\dots,K_c$ 
\EndProcedure
\end{algorithmic}
\end{algorithm}

{\bf Complexity:} Algorithm~\ref{alg:kmeans}: Searching the $K_n$-nearest neighbors of a data point using k-d trees takes $O(K_n\log J)$ in average and $O(K_nJ)$ in the worst case. Steps 3-5 has thus a complexity of $O(K_nJ\log J)$ in average and $O(K_nJ^2)$ in the worst case. Initial labelization is a simple linear partitioning of the 2D space in $K_c$ zones and is thus performed in $O(J)$. In the main loop, there are at most $O(MJK_c)$ calls to the function \texttt{FIT}. The Levenberg-Marquardt Algorithm requires $O(\epsilon^{-2})$ iterations to reach an $\epsilon$-approximation of a stationary point of the objective function~\cite{ueda2010global,bierlaire.15.book}. The overall complexity of Algorithm~\ref{alg:kmeans} is thus $O(K_nJ^2+MJK\epsilon^{-2})$.

\begin{figure*}[t!]
        \centering
        \includegraphics[width=0.6\linewidth]{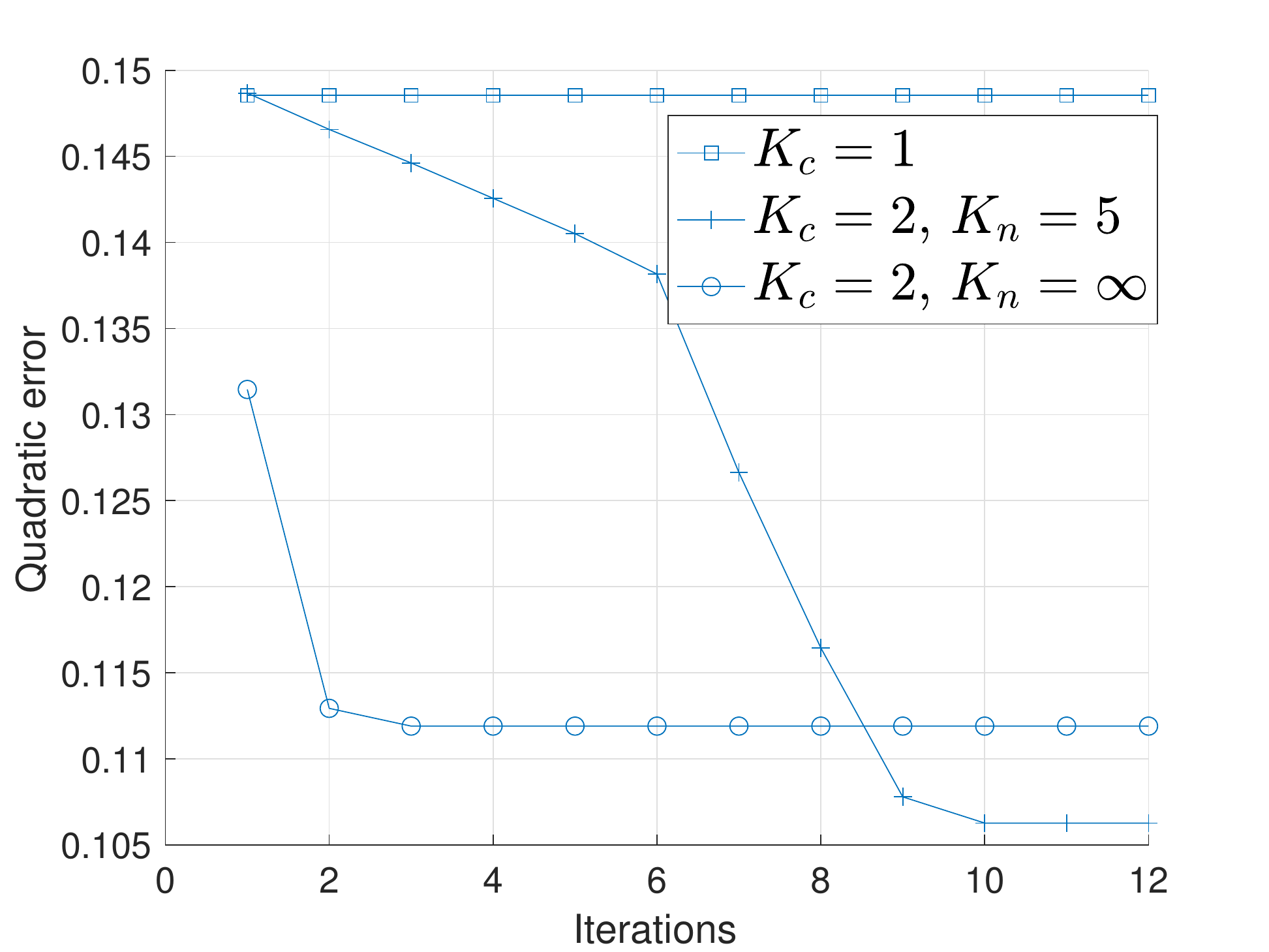}
    \caption{Quadratic error between normalized smoothed data and quadratic models.
    \label{fig:quaderror}
    }
\end{figure*}

\begin{figure*}[t!]
    \centering
    \begin{subfigure}[t]{0.5\textwidth}
        \centering
        \includegraphics[width=\linewidth]{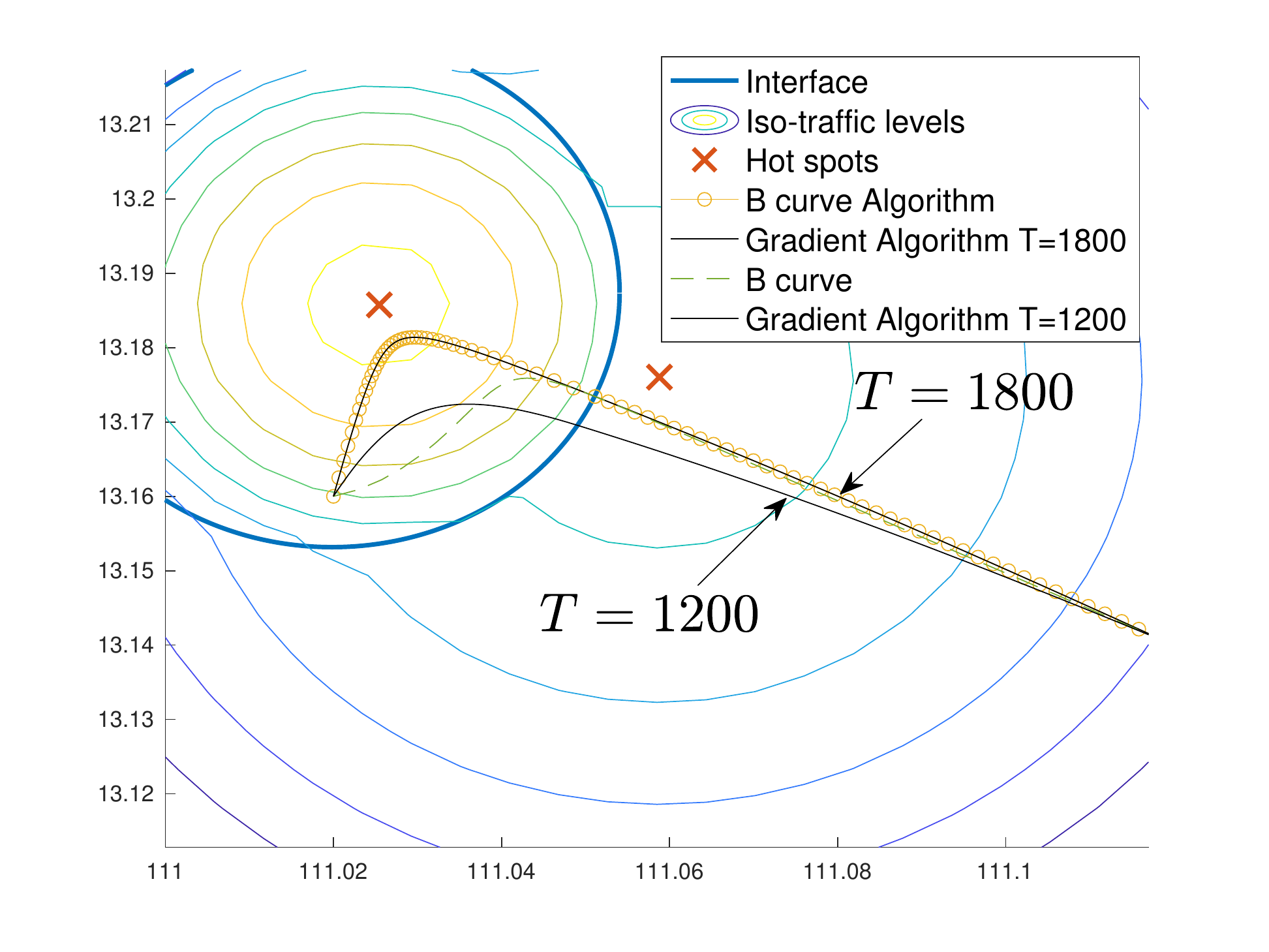}
        \caption{$K=60$, $T=1200$ and $1800$, 2 phases.}
        \label{fig:influenceofT2phases}
    \end{subfigure}%
    ~ 
    \begin{subfigure}[t]{0.5\textwidth}
        \centering
        \includegraphics[width=\linewidth]{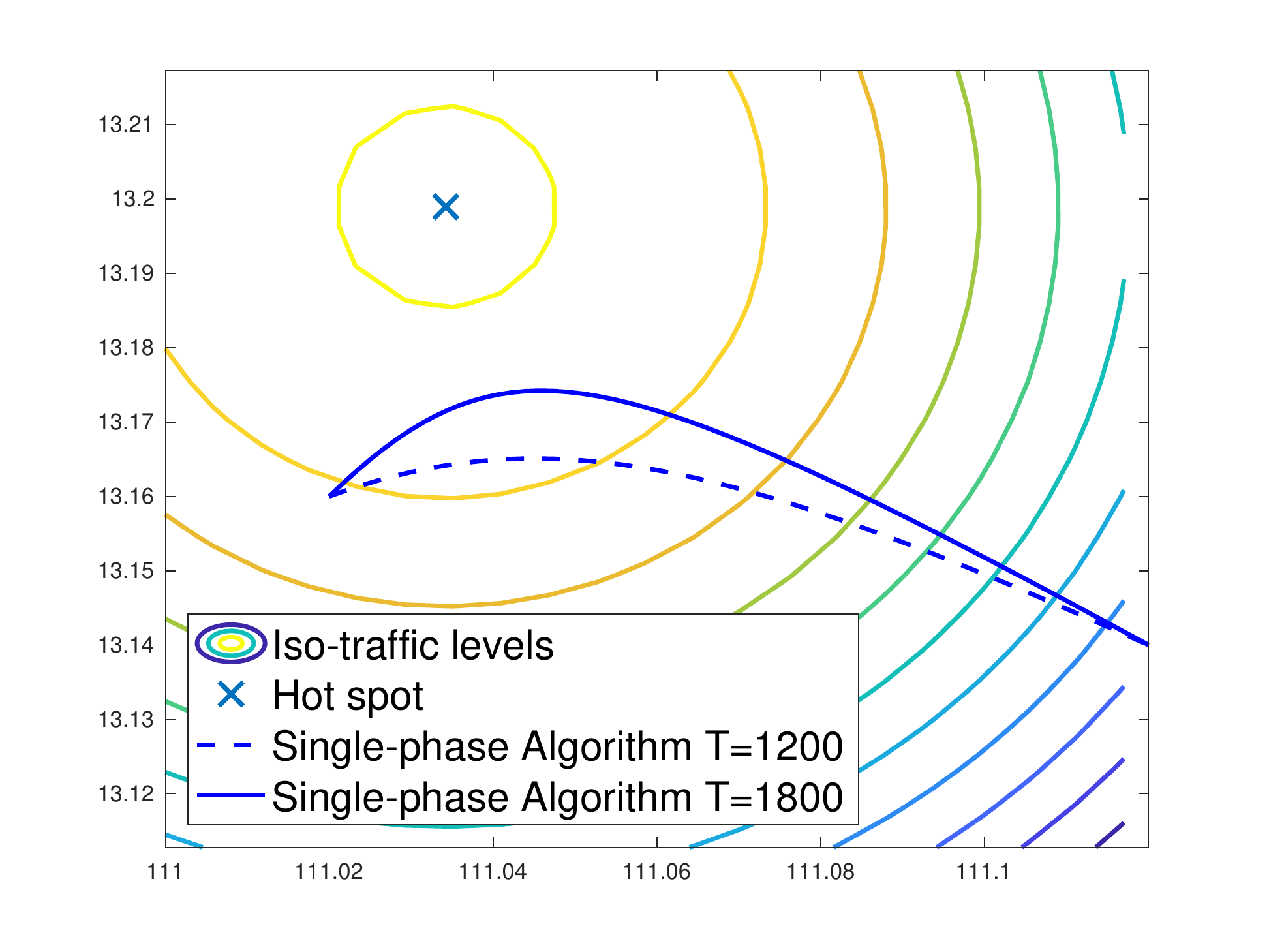}
        \caption{$K=60$, $T=1200$~s and $1800$~s, 1 phase.}
        \label{fig:influenceofT1phase}
    \end{subfigure}
    \caption{Influence of $T$ on the optimal trajectory.}
    \label{fig:influenceofT}
\end{figure*}

The K-means partition obtained after $M=12$ iterations is shown in Figure~\ref{fig:kmeanspartition}. The final fits are shown in Figures~\ref{fig:influenceofT2phases} and \ref{fig:influenceofT1phase} for $(K_c,\alpha)=(2,0.25)$ and $(K_c,\alpha)=(1,0.5)$, respectively. Figure~\ref{fig:quaderror} shows the quadratic error as a function of the number of iterations of the K-means algorithm for $K_c=1$ and $K_c=2$. The error is constant for $K_c=1$ as there is only one iteration, which performs the non-linear least square fitting for the single cluster. For $K_c=2$, we distinguish two cases: $K_n=5$ and $K_n=\infty$. In the former case, only the $5$ nearest neighbors of a data point are inspected to decide if a relabelization should be performed. In the later case, relabelization is systematically considered. Increasing $K_n$ increases the complexity of every iteration but provides a faster convergence.   

\begin{figure*}[t!]
    \centering
    \begin{subfigure}[t]{0.5\textwidth}
        \centering
        \includegraphics[width=\linewidth]{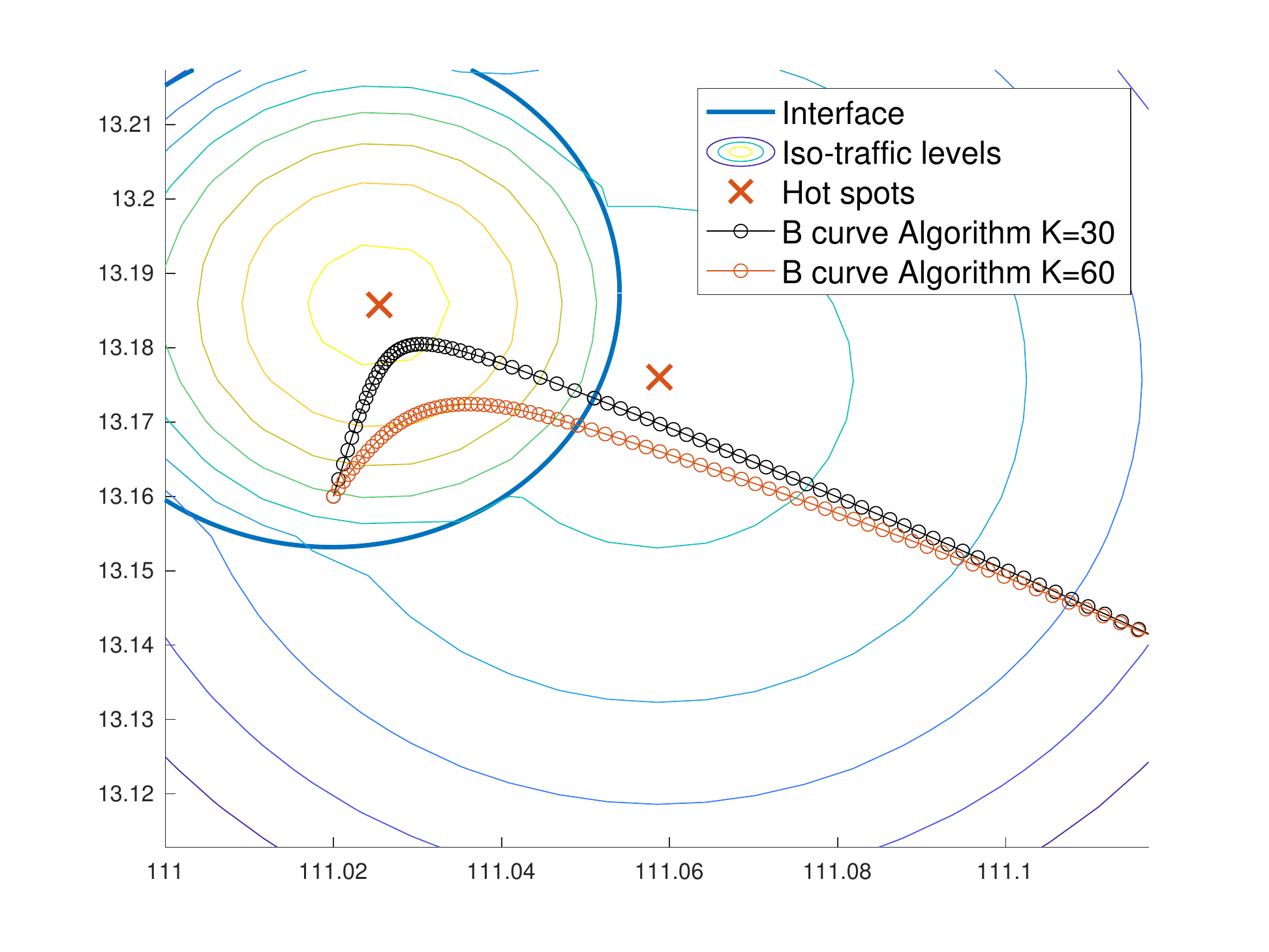}
        \caption{$K=30$ and $60$, $T=1200$, 2 phases.}
        \label{fig:influenceofK2phases}
    \end{subfigure}%
    \begin{subfigure}[t]{0.5\textwidth}
        \centering
        \includegraphics[width=\linewidth]{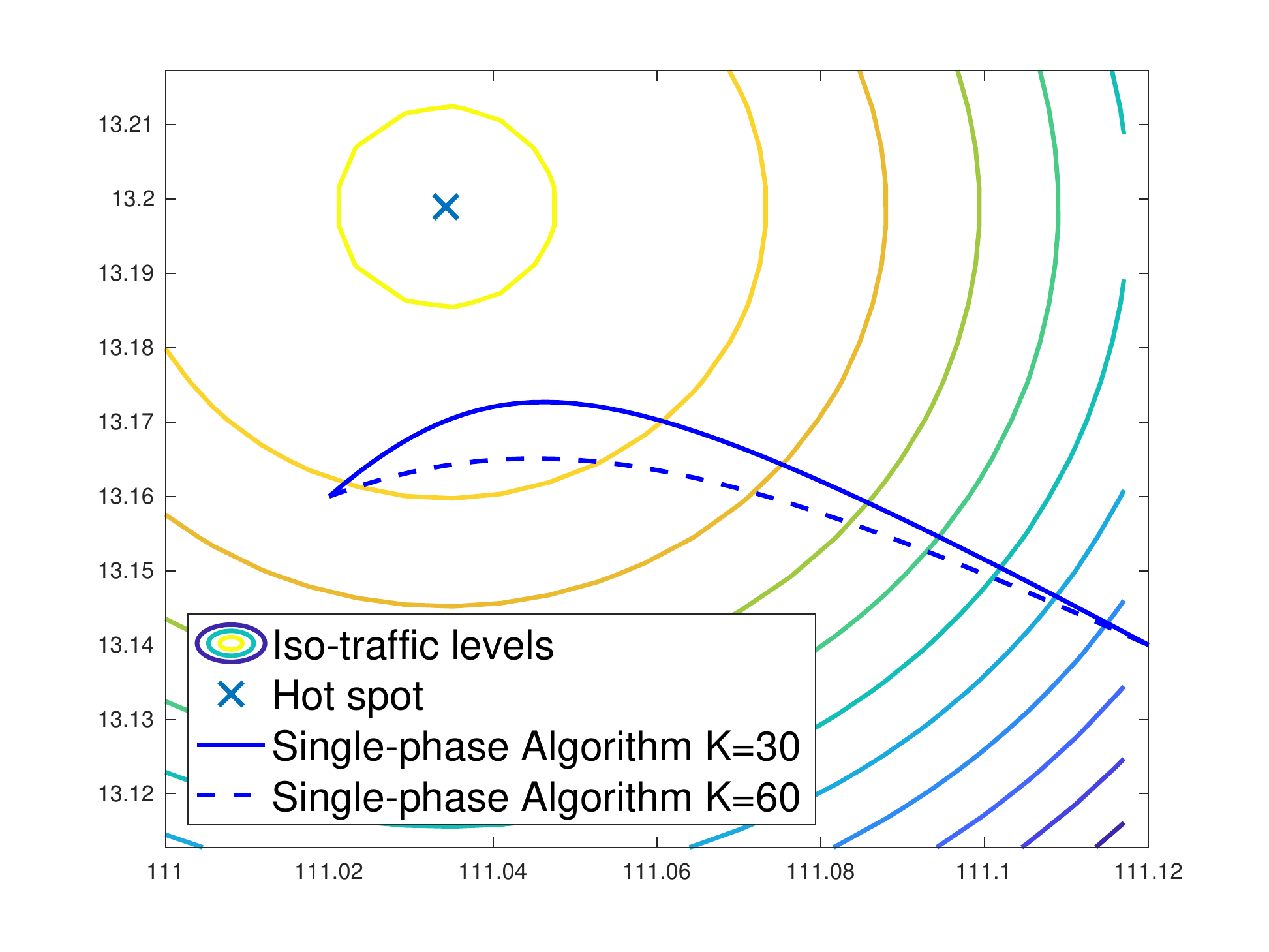}
        \caption{$K=30$ and $60$, $T=1200$~s, 1 phase.}
        \label{fig:influenceofK1phase}
    \end{subfigure}
    \caption{Influence of $K$ on the optimal trajectory.}
    \label{fig:influenceofK}
\end{figure*}
\subsection{Trajectory optimization results
\label{subsec:marceau-simus}}
\def\myvec#1#2{({#1}, {#2})}
\def\MYSHIFT{\myvec{111.}{13.}}
\subsubsection{Estimation procedure of
parameters K and T} 
as in most similar algorithms,
a 'prior' estimation of parameters
is necessary 
since the optimal results strongly depend on them.
Traffic parameters have been estimated 
just above,
so the required 
parameters are first the mass $K$
representing trade-off between the kinetic term 
and the traffic term in Lagrangian  
\eqref{eq:cost}
and then
the total duration time $T$ 
which has also a 
significant influence on the optimal trajectory
(by convention $t_0 = 0$).
This procedure
is developed in Appendix
\ref{estim-K-and-T:app} 
and yields in our case: 
\[
T \approx 850 
,
K 
= 50
\times r^2
\mbox{ with }
r = \frac{40 \times 10^6}{360 \times 100}
\mbox{ and we select: }
T = 1200, 1800s
\mbox{ and }
K =  30, 60 \times 
r^2.
\]

\subsubsection{Results}
\input{marceau-simus-6.tex}
optimal trajectories
are shown
in Figs. 
\ref{fig:influenceofT} 
and
\ref{fig:influenceofK}
and related numerical results  
in Table \ref{simus-marceau6:table}.
\\
Length of various trajectories
and related average 
velocities
are given 
in Table \ref{tab:length-speed}.
\\
The following main comments 
can be made on these results:
\begin{itemize}
\item First the obtained trajectories 
for $\GRADALGO$ and $\BALGO$
satisfy  both the optimal conditions, 
as
``relative" discrepancies
between impulsions and between Hamiltonians
at the interface are indeed below 
selected precision  
of  $2. 10^{-4}$,
numerically yielding a  full trajectory,
\item Then \GRADALGO
is found to be stable due to our choice of initial and terminal drone positions
    (not too close to the interface and not too far from their related hotspots) and of time intervals (not too large {\em temporal} phases). 
    Thus, we are far
    from the non-convexity conditions expressed 
    in Theorems
    \ref{th:single-phase-hessian}
    and
    \ref{th:two-phase-convexity}.
    The positivity of the total hessian is indeed constantly checked 
    at each iteration.
    Hence, 
    a large (spatio-temporal) convergence basin 
    resulting into exact convergence of the \GRADALGO (as well as the \BALGO)
    towards
    the {\em unique} optimal solution.
    \item 
    As already noted for the single-phase case 
    (paragraph 
    \ref{subsub:hotspot-hole}),
    decreasing  
    mass $K$  
    enables the drone to collect more traffic.
    In fact,
    during the allowed time interval, 
    the drone  will get closer to the hotspot with {\em maximal} traffic
    in order to decrease the total value function,
    \ie zone $1$ here.  
\item 
Increasing time interval 
($T\nearrow$) 
will produce the same tendancy,
allowing the drone to spend more time near
the hotspot with higher traffic 
$z_{h1}$, see Figs.
\ref{fig:influenceofT2phases}
and
\ref{fig:influenceofK2phases}.
\item Since
the drone has enough time to pick up traffic located close to hotspots, this 
is reflected in an average drone velocity
about half its nominal value
(see Table \ref{tab:length-speed}).
\end{itemize}
%

\section{Conclusion} \label{sec:conclusion}
In this paper, we propose a Lagrangian approach to solve the UAV base station optimal trajectory problem. When the traffic intensity exhibits a single phase, closed-form expressions for the trajectory and speed are derived from Hamilton-Jacobi equations. When the traffic intensity exhibits multiple phases, we characterize the crossing time and location at the interface. We propose two low-complexity algorithms for the bi-phase time-stationary traffic case that provide optimal crossing time and location on the interface and fulfill the necessary conditions of optimality. At last, we present a data processing procedure based on a modified K-means algorithm that derives a single-phase or bi-phase quadratic model from real traffic data. Further extensions of this work are envisioned to generalize the approach to three or more hot-spots and to consider multi-drone coordinated trajectories. 

\appendix
\subsection{Proof of Lemma~\ref{lemma:EulerLagrange}} \label{app:EulerLagrange}

In a neighborhood of the optimal trajectory, the first order variation of $S$ is null
\begin{eqnarray}
\delta S &=&\int_{t_0}^{T}\delta \mathcal{L}(t,z,a)dt \\ \notag
&=& \int_{t_0}^{T} \left[\nabla_z\mathcal{L}(t,z,a)\cdot\delta z(t) + \nabla_a\mathcal{L}(t,z,a)\cdot \delta a(t)  \right] dt. \\ \notag
\end{eqnarray}
We now note that $\delta a=\delta \frac{dz}{dt}=\frac{d(\delta z)}{dt}$. Integrating by part the second term in the integral of $\delta S$, we obtain
\begin{eqnarray}
\lefteqn{\int_{t_0}^T \nabla_a\mathcal{L}(t,z,a) \cdot \frac{d(\delta z)}{dt}dt}  \\ &=& [\delta z(t) \cdot \nabla_a\mathcal{L}(t,z,a)]_{t_0}^T - \int_{t_0}^T \delta z(t) \cdot\frac{d}{dt}\nabla_a \mathcal{L}(t,z,a)dt. \notag
\end{eqnarray}
Note that $[\delta z  \frac{\partial \mathcal{L}}{\partial a}]_{t_0}^T=0$ because $z_0$ and $z_T$ are fixed. 
Equating $\delta S$ to zero gives 
\begin{eqnarray}
0=\int_{t_0}^T \left[\nabla_z\mathcal{L}(t,z,a)-\frac{d}{dt}\nabla_a\mathcal{L}(t,z,a)  \right]\cdot \delta z(t) dt.
\end{eqnarray}
As this should be true for every $\delta z$, $\mathcal{L}$, $z_0$ and $z_T$, we obtain the first result. 

Now assume that we have the optimal $a(t)$, the condition for $z(T)$ to be the optimal final position is
\begin{eqnarray}
\delta S & =& [\delta z(t)\cdot  \nabla_a \mathcal{L}(t,z,a)]_{t_0}^T + \nabla J(z(T))\cdot\delta z(T) \notag \\
&=& \nabla_a\mathcal{L}(z(T),T,a(T))\cdot \delta z(T) + \nabla J(z(T))\cdot \delta z(T) \notag \\
&=&0.
\end{eqnarray}
Note that $z_0$ is fixed and so $\delta z$ in $z_0$ is null. We thus obtain the second result of the lemma. 

\subsection{Proof of Lemma~\ref{lemma:lagrangienhomogene}} \label{app:lagrangienhomogene}

As $\mathcal{L}(z,a)$ is an homogeneous function of $z$ and $a$, we have: $\mathcal{L}(\lambda z,\lambda a)=|\lambda|^{\alpha}\mathcal{L}(z,a)$ for all $\lambda$ (in our case with $\alpha=2$). Deriving this expression with respect to $\lambda$, setting $\lambda=1$, and noting that $a=\dot{z}$ we obtain
\begin{eqnarray}
z\cdot \frac{\partial \mathcal{L}(z,\dot{z})}{\partial z} +\dot{z}\cdot \frac{\partial \mathcal{L}(z,\dot{z})}{\partial z}&=&\alpha \mathcal{L}(z,\dot{z}). \label{eq:eulerhom}
\end{eqnarray}
Using (\ref{eq:eulerimpulsion}) and (\ref{eq:eulerhom}), we have:  $z\cdot \frac{d p}{dt}+\dot{z}\cdot p=\alpha \mathcal{L}$ or equivalently $\frac{d(p\cdot z)}{dt}=\alpha \mathcal{L}$. We can now integrate the cost function (\ref{eq:problem}) along the optimal trajectory as follows
\begin{eqnarray}
S(t_0, z_0,T,z_T)&=&\frac{1}{\alpha} \int_{t_0}^T \frac{d(p\cdot z)}{dt}(t) dt+J(z_T) \notag \\
&=&\frac{1}{\alpha} \left(p(T)\cdot z_T-p(t_0)\cdot z_0\right ) +J(z_T) 
\label{eq:homogeneous-result-S}
\end{eqnarray}

\subsection{Proof of Lemma~\ref{lemma:hamiltonjacobi}}
\label{app:hamiltonjacobi}

We assume that an optimal trajectory exists and we apply the principle of optimality on the optimal trajectory between $(t,z^*(t))$ and $(t+h,z^*(t)+ah)$, where $h>0$, to opbtain
\begin{eqnarray}
S(t,z^*(t),
T,z_T)  
&=& \min_{a} [h \mathcal{L}(z^*(t),a)+S(t+h,z^*(t)+ah,
T,z_T)] \notag \\ 
&=& \min_{a}[h \mathcal{L}(z,a)+ \notag
S(t,z^*(t),T,z_T) 
+ha\cdot \nabla_{\Xone} S(t,z^*(t),T,z_T)+ \notag 
h\frac{\partial S}{\partial \Tone}(t,z^*(t),
T,z_T)]. \notag 
\end{eqnarray} 
This implies that
\begin{eqnarray}
\frac{\partial S}{\partial 
\Tone}(t,z^*(t),T,z_T) 
&=& -\min_{a} [a\cdot \nabla_{\Xone} S(t,z^*(t),T,z_T)+\mathcal{L}(z^*(t),a)] \notag   \\
&=& \max_a [-a\cdot \nabla_{\Xone} S(t,z^*(t),T,z_T)-\mathcal{L}(z^*(t),a)] \notag \\
&=& H(t,z^*(t),
-\nabla_{\Xone} 
S(t,z^*(t),T,z_T)). \notag 
\end{eqnarray}
By using the same approach between $t-h$ and $t$, we deduce in the same way equation \eqref{Hamilton-Jacobi-forward-ms:eq} when the final time $T$ is varying.  

\subsection{Proof of Theorem~\ref{th:quadraticfunction}} \label{app:quadraticfunction}

From (\ref{eq:eulerlagrange2}) and (\ref{eq:Lsinglephase}), we obtain the following ordinary differential equation of second degree: $\ddot{z}=-\frac{u_0}{K}z$. If $\frac{u_0}{K}>0$, we define $\omega^2=\frac{u_0}{K}$ and we look for an optimal trajectory of the form: $z(t)=A\cos(\omega t)+B\sin(\omega t)$. If $\frac{u_0}{K}<0$, we look for an optimal trajectory of the form: $z(t)=A\cosh(\omega t)+B\sinh(\omega t)$ with $\omega^2=-\frac{u_0}{K}$. Let us denote $z_0=z(t_0)$ and $a_0=a(t_0)$ the initial conditions for $z$ and $\dot{z}$. 

Take the case $\frac{u_0}{K}<0$. Using the derivative of $z(t)$ and identifying terms, we obtain: $z(t)=z_0\cosh \omega(t-t_0) + \frac{a_0}{\omega}\sinh \omega(t-t_0)$. At $t=T$, we have also: $z_T=z_0\cosh \omega(T-t_0)+\frac{a_0}{\omega}\sinh \omega(T-t_0)$, from which we deduce
\begin{eqnarray}
a(t_0)&=&\frac{\omega(z_T-z_0\cosh \omega(T-t_0))}{\sinh \omega(T-t_0))}, 
\label{marc-ptzero:eq}
\\
a(T)&=&\frac{\omega(-z_0+z_T\cosh \omega(T-t_0))}{\sinh \omega(T-t_0))}.
\label{marc-pT:eq}
\end{eqnarray}
when $u_0<0$. In a similar way, we have
\begin{eqnarray}
a(t_0)&=&\frac{\omega(z_T-z_0\cos \omega(T-t_0))}{\sin \omega(T-t_0)}, \\
a(T)&=&\frac{\omega(-z_0+z_T\cos \omega(T-t_0))}{\sin \omega(T-t_0)},
\end{eqnarray}
when $u_0>0$. Injecting $a(t_0)=a_0$ in the equation of the trajectory provides the result.

For the computation of $S$, we now use the result of Lemma~\ref{lemma:lagrangienhomogene} as our cost function is 2-homogeneous.
From equation (\ref{eq:valuefunction}), we see that only initial and final conditions are required to compute the cost function. Recall now that $p=Ka$. 
Injecting the equations of $a(t_0)$ and $a(T)$ in (\ref{eq:valuefunction}), we obtain the result for the cost function. 
%
%

\subsection{Proof of Theorem~\ref{th:multiphase}}
\label{app:multiphase}

We assume that the location and time $(\xi,\tau)$ of interface crossing is known and unique. The optimal trajectory between $(z_0,t_0)$ and $(z_T,T)$ can be decomposed in two sub-trajectories that are themselves optimal between $(z_0,t_0)$ and $(\xi,\tau)$ on the one hand and between $(\xi,\tau)$ and $(z_T,T)$ on the other hand, by the principle of optimality. 
In region 1, the optimal cost up to $\tau$ is
\begin{equation}
S_1(t_0,z_0,\tau,\xi)=\int_{t_0}^{\tau}\mathcal{L}(z^*(s),a^*(s))ds.
\end{equation}
Using Hamilton-Jacobi, we obtain
\begin{equation}
\frac{\partial S_1}{\partial T_2} 
(t_0,z_0,\tau,\xi)=-H_1(\xi,p^*(\tau^-)). 
\end{equation}
In the same way, the optimal cost in region 2 is
\begin{eqnarray}
S_2(\tau,\xi,T,z_T)&=&\int_{\tau}^{T}\mathcal{L}(z^*(s),a^*(s))ds. 
\end{eqnarray}
Using again Hamilton-Jacobi, we obtain
\begin{equation}
\frac{\partial S_2}{\partial T_1} 
(\tau,\xi,T,z_T)=H_2(\xi,p^*(\tau^+)).
\end{equation}
A necessary condition for the optimality of $\tau$ is thus
\begin{equation}
\frac{\partial S_1}{\partial \Ttwo}(t_0,z_0,\tau,\xi)+\frac{\partial S_2}{\partial \Tone}(\tau,\xi,T,z_T)=0,
\label{optimal-tau-first-ms:eq}
\end{equation}
that is
\begin{equation}
H_1(\xi,p^*(\tau^-))=H_2(\xi,p^*(\tau^+))\label{optimal-tau-second-ms:eq}.
\end{equation}
A necessary condition for the optimality of $\xi$ in 
the total cost
under the constraint $f(\xi)=C$ 
is also 
\begin{eqnarray}
\mu\ \nabla_z f(\xi)
&=&
\nabla_{\Xtwo} 
S_1(t_0,z_0,\tau,\xi)
+
\nabla_{\Xone}
S_2(\tau,\xi,T,z_T) \notag
\\ &=&
p^*(\tau^-) - p^*(\tau^+) \notag
\end{eqnarray}
where $\mu$ is a Lagrange multiplier associated to the constraint and where
the second line comes from
equation (\ref{eq:optimpulsion})
of Hamilton-Jacobi. 
Thus we obtain precisely equation (\ref{impulsion-ms:eq}).

\subsection{The structure and the diagonalization
of single- and two- phase Hessian of the value fonction}
\label{app:Hessians}
\subsubsection{
Proof of Theorem
\ref{th:single-phase-hessian}
\label{structure-hessian:app}}
Recall that 
we study the second-order differentiability properties of
single-phase value function
$\myS$.
Recall that from Eq. 
\eqref{eq:definition-hessian}
we study Hessians of the form
\def\abbrevS{\sigma}
\[
\mathcal{M}_{1+2, 1+2}(\bR) \ni \nabla^2
\abbrevS
=
\begin{pmatrix}
\ddpart{}{T_i}\
\abbrevS
& 
\dpart{}{T_i}\nabla_{Xi}\ \abbrevS 
\\
\nabla_{Xi}\dpart{}{T_i}\
\abbrevS 
& \nabla^2_{X_i,X_i}\
\abbrevS
\end{pmatrix}
\ \
\mbox{where } \abbrevS 
\mbox{ stands for } \myS
\mbox{ and }
i = 1,2
\]
It is clear by inspecting 
the symmetries of (\ref{eq:vu0neg}) 
with respect to the spatial coordinates 
that  
\[
\nabla^2_{\Xone,\Xone} \myS
=
\nabla^2_{\Xtwo,\Xtwo} \myS
=
K g\ \Id_2 
,
\]
\begin{equation}
\mbox{where }
g = \w\ \coth 
\phase
\ \ \ 
\mbox{and}
\ \phase = \w\ (\ttwo - \tone) 
\mbox{ is the {\em temporal} phase}
.
\label{temporal-phase:eq}
\end{equation}
Concerning symmetry with respect to time variables, 
one also finds easily that
\[
\ddpart{S}{\Tone}(\tone, \xone, \ttwo, \xtwo) 
= 
\ddpart{S}{\Ttwo}(\tone, \xone, \ttwo, \xtwo). 
\]
The Hessian 
with respect to the variable 
$\chi_i$ = 
$(T_i, X_i)_{i = 1,2} \in  \bR \times \bR^2$ 
has thus the following structure:
\begin{align}
{\cal H}(\chi_i) 
&=
\Hessian{\alpha}{K\ g}{\Pi_i}
\hspace*{15mm}
\mbox{with}
\hspace*{15mm}
\left\{
\begin{array}{l}
\alpha = \ddpart{S}{\Tone}
= \ddpart{S}{\Ttwo}
\\[4mm]
g = \w\ \coth \phase
\\[4mm]
\Pi_i = 
\dapart{\nabla_{X_i} S}{T_i} 
\end{array}
\right.
\label{Hessian-first:eq}
\end{align}
\\[1mm]
where the vector $\Pi_i$
and  real scalar $\alpha$ 
remain to be determined.
\\
Now
recall that initial and final impulsions 
$\pone$  
and $\ptwo$ 
as given by
\eqref{marc-ptzero:eq}
write: 
\begin{align}
\pone&= K\ \w\ \dfrac{-\xone \cosh \phase+\xtwo}{\sinh \phase}, 
\label{pone:eq}
\\
\ptwo&= K\ \w\ \dfrac{-\xone+\xtwo \cosh \phase}{\sinh \phase}.
\label{ptwo:eq}
\end{align}
In fact the expression
of $\Pi_i$ and $\alpha$ 
appearing in 
Hessian
\eqref{Hessian-first:eq}
result 
both 
from the computation of
$\dpart{p_i}{T_i}$   
for $i = 1,\ 2.$
Using \eqref{temporal-phase:eq},
the differentiation of
\eqref{pone:eq} 
and \eqref{ptwo:eq}
\wrt $\Tone$ 
\resp $\Ttwo$ 
leads to the following simple result
\begin{equation}
\label{dp/dt:eq} 
\Pi_1 =
\dapart{\nabla_{\Xone} S}{\Tone} =
\dpart{\pone}{\Tone} = 
-\ \w\ \dfrac{\ptwo}{\sinh \phase}
 \ \ \ 
 \mbox{ resp. }
 \ \ 
 \Pi_2 =
\dapart{\nabla_{\Xtwo} S}{\Ttwo} =
 \dpart{\ptwo}{\Ttwo} =
 -\ 
 \w\ \dfrac{\pone}{\sinh \phase}.
\end{equation}
Now, as far as the
second partial derivative of value function
$S$ 
with respect to 
time ($\alpha$) is concerned 
\begin{equation}
\alpha = \ddpart{S}{\Tone} = - \dpart{H}{\Tone}
\ \ \ 
\mbox{ (Hamilton-Jacobi) }
\label{eq:HJ-time}
\end{equation}
Recalling that the Hamiltonian for any Newtonian model 
has the form
\begin{equation}
H(z, p, t)  = 
\dfrac{{\norm{p}}^2}{2\ K} + V(z)
\label{eq:newton-hamiltonian}
\end{equation}
and differentiating it with respect to the 
{\em  temporal phase} 
$\phase$
at {\em fixed} extremities $\xone$ and $\xtwo$
leads to
\begin{align} 
\del H 
&= \dfrac{\pone \cdot \del \pone}{K}
=  \dfrac{\ptwo \cdot \del \ptwo }{K}
= - \dfrac{\pone \cdot \ptwo} {K\ \sinh \phase}\ \del \phase
\enspace,
\ \
\mbox{that is}
\nonumber
\\[2mm]
\alpha &= \ddpart{S}{\Tone} 
= \ddpart{S}{\Ttwo} 
=  
-\ \w \ 
\dpart{H}{\phase} = \dfrac{\w}{K}\ 
\dfrac{\pone \cdot \ptwo} {\sinh \phase}.
\qed
\label{d2S-dphase2:eq}
\end{align} 
Now 
from \eqref{eq:sumcost}
the two-phase Hessian is simply the sum of 
the two single-phase Hessians.

\subsubsection{Diagonalizing the Hessian of 
single- and two-
phase value function 
\label{diagonal-hessian:app}}
since the Hessians of single- and two-phase value function do have the same structure,
we address a 
general
Hessian of the form
\eqref{Hessian-first:eq}.
\\
Its characteristic polynom  
is easily developed as 
\[
\det\ ({\cal H} - \nu\ \Id_{3})
=
(K g - \nu)^2\ (\alpha - \nu)  
-\, (K g - \nu)\ \norm{\Pi}{}^2.
\]
Since
from \eqref{Hessian-first:eq}
one has
$
\forall 
\trivec{\theta} {x}\in \bR \times \bR^2\
\ \ 
{\cal H} \trivec{\theta}{x} 
= 
\trivec{\alpha\ \theta + \Pi \cdot x}
{\Pi\ \theta + K g\ x}
$, 
\\
the three real eigenvalues 
$(\nu_i)_{i = 0:2}$
and related eigenvectors 
$Y_i = (\theta_i,x_i)$
thus satisfy
\\[1mm]
$
i)\ \text{if }\nu_0 = K g > 0 \text{ then}
\\
Y_0 =
\trivec{\theta = 0}{x_0} 
\mbox{ with } x_0 \perp \Pi
\enspace:
\mbox{"pure space" eigenvector }
Y_0
\in \{0\} \times \bR^2 ,
\\[1mm]
\hspace*{-2mm}
ii)\ \text{if }
\nu_1,\ \nu_2
\mbox{ such that }
(K g - \nu)\ (\alpha - \nu)\ - \norm{\Pi}{}^2 
= 0 \text{ then }
\\
%
Y_i = {\trivec{\theta_i}{x_i = \Pi}}_{i = 1,2} 
\enspace:
\mbox{two 
(mutually orthogonal)
eigenvectors orthogonal to } 
Y_0 
\ (\mbox{since }  \Pi \perp x_0 ) 
\\
\mbox{  with  } 
\nu_1\ \nu_2 
= 
\alpha\ K g\ - \norm{\Pi}{}^2
 \ \
= \det({\cal H})\ /\ (K g)
\ \ 
\mbox{ and }
\ \ 
\nu_1 + \nu_2 = \alpha + K g.
$
\\[1mm]
Thus
$\nu_1, \nu_2$ cannot be simultaneously
negative  since this would imply
\\
$\alpha < - K g < 0 \Ra \det({\cal H})  = \nu_0\ \nu_1\ \nu_2\ < 0$.
\\
However the condition
$\alpha\ K g\  <\ \norm{\Pi}{}^2$ 
(\mbox{e.g.} induced by the {\em sufficient}
condition 
$\alpha < 0$) implies that
one of eigenvalues 
$\nu_1\ ,\ \nu_2 < 0$
\ \ 
namely the {\em local non-convexity}
of the Hessian matrix
${\cal H}$.
\\
This property holds thus for each single-phase
value function as well as for the total
two-phase case.

\subsection{Estimation of parameters}
\label{estim-K-and-T:app}
It is based on the following remarks.
First, the interface location ${\cal I}$ 
is independent of
$K$ and $T$ 
(since only depends on {\em traffic} variables).
Then, final time $T$ 
remains always estimated in \underline{seconds}
$(t_0 = 0)$.
Also, according to \cite{fotouhi2017dronecells}, a typical drone cell flights at 
speed 
$\bar{V} = 20$~m/s 
and has an autonomy of about $28$~min.
Last, we test our procedure
for spatially-scaled data, 
where the convenient
spatial unity is 100\ m,
so that the scaling ratio is
\begin{equation}
r = \frac{40 . 10^6}{360} \times 100 
\approx 1111.1111 
\label{eq:r-scaling}
\end{equation}
Also, the following numerical estimations are considered
\renewcommand{\labelitemi}{\theenumiii}
\begin{itemize}
\item[a)] 
for estimating $T$:
the total length of the trajectory
$L$ is approximated by
\[
L 
\approx 
c\ \times \parallel z_0  - z_T \parallel
\mbox{ with }
c = 1.5
\ra T = \frac{\bar{V}} {L}
\]
(recall that $T$ is un-scaled
and stands in seconds.)
\item[b)] 
for estimating $K$:
the phase 
in each sub-trajectory should satisfy
\\
$\phase_i = \omega_i\ T_i<
\phase_{max}
\enspace, \ 
i = 1,2$
where 
$\phase_{max} 
=
10$
(
also an un-scaled constant)
.
\\
Then, spatial scaling by 
the quantity $r$ 
implies
\[
u_0^i \gets 
\dfrac{u_0^i}{r^2}
\enspace,
i = 1,2
\enspace, 
\mbox{ and thus: }
K \gets \dfrac{K}{r^2}
\mbox{ since }
\ \ 
K\ \omega_i ^2 = u_0^i
\ 
\mbox{ at any scale.}
\]
\end{itemize}
Passing back into the original
frame therefore implies that 
$ K \gets K\ r^2$. 



\bibliographystyle{IEEEtran}
\bibliography{IEEEabrv,dronebs}

\end{document}

%% file: marc-macros.tex
\def\wrt{\mbox{wrt. \,}}
\def\resp{\mbox{ resp. }}

\def\ra{\rightarrow}
\def\Ra{\Rightarrow}

\def\Xone{X_1}
\def\Tone{T_1}
\def\Xtwo{X_2}
\def\Ttwo{T_2}

\def\Pone{\nabla_{\Xone}}

\def\Hone#1{\dpart{#1}{\Tone}} 
\def\Htwo#1{\dpart{#1}{\Ttwo}} 

\def\xx{x} 			
\def\indexone{1} 	
\def\indextwo{2}	 

\def\xone{{\xx}_\indexone}
\def\xtwo{{\xx}_\indextwo}

\def\pp{p}
\def\pone{{\pp}_\indexone}
\def\ptwo{{\pp}_\indextwo}

\def\tt{t}
\def\tone{\tt_\indexone} 
\def\ttwo{\tt_\indextwo} 

\def\w{\omega}
\def\phase{\Phi}
\def\phase{\phi}

\def\transp#1{{#1}^{\dagger}}

\def\norm#1{\parallel {#1} \parallel}

\def\trivec#1#2{
\Big(
\begin{array}{c}
{
{#1}} 
\\[-2mm] 
\rule{2pt}{0.1pt} 
\\[-1.5mm] 
{
{#2}} 
\end{array} \Big)}

\def\trivec#1#2{({#1},\,{#2})} 

\def\dpart#1#2{\dfrac{\partial {#1}}{\partial {#2}}}
\def\dapart#1#2{\dfrac{\partial}{\partial #2} {#1}}
\def\ddpart#1#2{\dfrac{\partial^2 {#1}}{\partial {#2}^2}}

\def\del{\delta\,}

\def\myS{S(\tone, \xone, \ttwo, \xtwo)} 
\def\Id{\mbox{Id}}
\def\Real{\bR} 

\def\GRADALGO{\mbox {\textsc{Grad-Algo} }}
\def\BALGO{\mbox {\textsc{B-Algo}}}

\def\HSPONE{z_{h1}}
\def\HSPTWO{z_{h2}}

\def\Hessian#1#2#3{
\left(
\begin{array}{l | c}
#1
&
\hspace*{2mm} 
\transp{#3} 
\\[1mm]
\hline
#3 
&
\hspace*{4mm}
#2\ \Id_2 
\end{array}
\right)
}

\def\Vmean{\bar{V}}

%% file: marceau-simus-6.tex
\begin{table*}[t!] 
\[
\hspace*{-8mm} 
\begin{array}{|c 
c c c c 
| l 
| 
} 
\hline
\hline
\vspace*{1mm}
\xi
& {\mathbf S}  & H^-  
& \phi^- & \phi^+  
& \tau  
%
\\
\hline
\hline
\multicolumn{6}{|c|}
{ \mathbf{a)}\ 
K = 60 \times r^2
\enspace,
\ T = 1200
\ 
\ra S_{1ph} = -221.4103
}
\\[1mm]
\multicolumn{6}{|c|}{\BALGO\ (Niter = 12) }
\\
\myvec{111.048
}{13.170}
& -168.5713 
&  0.8131 
& 1.9803  & 0.9236 
&538.69
\\[1mm]
\hline
\multicolumn{6}{|c|}{\GRADALGO (Niter = 5)}
\\
\myvec{111.048
}{13.170}
& -168.5713  
& 0.8131 
& 1.9802 &  0.9237 
& 538.68
\\
\hline
\hline
\multicolumn{6}{|c|}
{ \mathbf{b)}\
K = 60 \times r^2
\enspace,
\ T = 1800
\
\ra S_{1ph} = -637.9077
}
\\[1mm]
\multicolumn{6}{|c|}{\BALGO\ (Niter = 12)}
\\
\myvec{111.0511}{13.1734}
& -641.6388  
&  0.8332 
& 4.1764  & 0.9273 
&1136.09
\\[1mm]
\hline
\multicolumn{6}{|c|}{\GRADALGO (Niter = 2)}
\\
\myvec{111.0511}{13.1734}
& -641.6388    
&  0.8332
& 4.1763 & 0.9273 
& 1136.08
\\
\hline
\hline
\multicolumn{6}{|c|}
{ \mathbf{c)}\
K = 30 \times r^2
\enspace,
\ T = 1200
\
\ra S_{1ph} = -399.1474 
}
\\[1mm]
\multicolumn{6}{|c|}{\BALGO\ (Niter = 12)}
\\
\myvec{111.051}{13.173} 
& -393.1217  
&  0.8313  
&  3.7936 & 0.9289 
& 729.71 
\\[1mm]
\hline
\multicolumn{6}{|c|}{\GRADALGO (Niter = 2)}
\\
\myvec{11.051}{13.173} 
& -393.1217   
&   0.8313 
&  3.7933 &  0.9290 
& 729.65 
\\
\hline
\end{array}
\]
\caption{\label{simus-marceau6:table}
Table of results with 
\BALGO\ and \GRADALGO 
for 
$K = 30, 60 \times r^2$
,
$T = 1200, 1800 s$
\\
and
$r = 
\frac{40 \times 10^6}{360 \times 100}
\approx 1111.1111$
(corresponding to scaled spatial unity of 100m).
\\
Also, $S_{1ph}$ denotes the total  single-phase cost
(\ie assuming only one, effective hotspot).
} 
\end{table*}
\begin{table*}[t!]
\[
\hspace*{-8mm} 
\begin{array}{|  l |  c c | c c  |}
\hline
\hline
\multicolumn{1}{| l |}
{case}
& 
\multicolumn{2}{ c |}
{\mbox{bi-phase}}
& 
\multicolumn{2}{ c |}
{\mbox{single-phase}}
\\[1mm]
&
L\ (\mbox{km}) & \Vmean (\mbox{m/s})
& 
L\ (\mbox{km}) & \Vmean (\mbox{m/s})
\\
\hline
 {\mathbf{a)}} 
 &
 12.41 & 10.34 
 &
 11.71 & 9.76 
\\[1mm]
\hline
{\mathbf{b)}} 
&
13.80 & 7.66 
 &
12.53 & 6.96
\\[1mm]
\hline
{\mathbf{c)}} 
&
13.64 & 11.37 
&
12.37 & 10.30 
\\
\hline
\end{array}
\]
\caption{\label{tab:length-speed} Table of trajectory lengths and related average speeds.}
\end{table*}